\documentclass[11pt]{amsart}
\usepackage{amssymb,euscript,tikz,units}

\usepackage[colorlinks, linkcolor=blue, anchorcolor=blue, citecolor=blue]{hyperref}

\usepackage{subcaption}

\usepackage{verbatim}
\usepackage[nameinlink]{cleveref}
\usepackage{colonequals}
\usepackage{tikz}
\usepackage{enumerate}
\usepackage[justification=centering]{caption}
\usepackage{amsfonts,amssymb,amsmath,amsthm,amscd}
\usepackage{bbm}

\allowdisplaybreaks

\numberwithin{equation}{section}

\newcommand{\cM}{\mathcal{M}}

\newcommand{\cC}{\mathcal{C}}

\newcommand{\cN}{\mathcal{N}}
\newcommand{\cO}{\mathcal{O}}
\newcommand{\cP}{\mathcal{P}}

\newcommand{\cS}{\mathcal{S}}
\newcommand{\cT}{\mathcal{T}}

\newcommand{\C}{\mathbb{C}}
\renewcommand{\H}{\mathbb{H}}
\newcommand{\N}{\mathbb{N}}
\renewcommand{\P}{\mathbb{P}}

\newcommand{\R}{\mathbb{R}}
\renewcommand{\S}{\mathbb{S}}
\newcommand{\T}{\mathbb{T}}
\newcommand{\W}{\mathbb{W}}

\newcommand{\eps}{\varepsilon}

\newcommand{\dvol}{d\mathrm{vol}}

\newcommand{\EgnL}{\mathbb{E}_{\rm WP}^{g,n,\mathbf{L}}}

\def\area{\mathop{\rm Area}}

\def\Mod{\mathop{\rm Mod}}

\theoremstyle{plain}
\newtheorem{theorem}{Theorem}[section]
\newtheorem{corollary}[theorem]{Corollary}
\newtheorem{proposition}[theorem]{Proposition}
\newtheorem{lemma}[theorem]{Lemma}

\newtheorem{remark}[theorem]{Remark}

\newtheorem*{thm*}{Theorem}

\theoremstyle{definition}

\theoremstyle{remark}

\newtheorem*{rem*}{Remark}
\newtheorem*{def*}{Definition}
\newtheorem*{con*}{Construction}
\newtheorem*{definition*}{Definition}

\newcommand{\vone}{\vskip 1\baselineskip}
\newcommand{\p}{\partial}
\newcommand{\bpmatrix}{\begin{pmatrix}}
\newcommand{\epmatrix}{\end{pmatrix}}
\newcommand{\supp}{\mathrm{Supp}}

\newcommand{\be}{\begin{equation}}
\newcommand{\ene}{\end{equation}}
\newcommand{\br}{\begin{remark}}
\newcommand{\er}{\end{remark}}
\newcommand{\bl}{\begin{lemma}}
\newcommand{\el}{\end{lemma}}
\newcommand{\bcor}{\begin{corollary}}
\newcommand{\ecor}{\end{corollary}}
\newcommand{\bpro}{\begin{proposition}}
\newcommand{\epro}{\end{proposition}}
\newcommand{\ben}{\begin{enumerate}}
\newcommand{\een}{\end{enumerate}}
\newcommand{\bp}{\begin{proof}}
\newcommand{\ep}{\end{proof}}
\newcommand{\beq}{\begin{equation*}}
\newcommand{\eeq}{\end{equation*}}
\newcommand{\bear}{\begin{eqnarray}}
\newcommand{\eear}{\end{eqnarray}}
\newcommand{\beqar}{\begin{eqnarray*}}
\newcommand{\eeqar}{\end{eqnarray*}}
\newcommand{\bt}{\begin{theorem}}
\newcommand{\et}{\end{theorem}}
\newcommand{\bex}{\begin{excer}}
\newcommand{\eex}{\end{excer}}

\makeatletter

\newcommand{\Rmnum}[1]{\expandafter\@slowromancap\romannumeral #1@}
\makeatother

\providecommand{\keywords}[1]{\textbf{Keywords:} #1}


\begin{document}

\title[The total mass of Brownian loop measure]{The total mass of Brownian loop measure of Riemann surfaces for large genus}
\author{Jiankun Hou and Yunhui Wu}

\vspace{.1in}

\address{Qiuzhen College, Tsinghua University, Beijing, China}
\email[(J.~H.)]{houjk21@mails.tsinghua.edu.cn}

\address{Yau Mathematical Sciences Center and Department of Mathematical Sciences, Tsinghua University, Beijing, China}
\email[(Y.~W.)]{yunhui\_wu@tsinghua.edu.cn}

%\address{}
\email[]{}
\email[]{}

\address{}
\email[]{}

\date{}


\maketitle
\vspace{-.2in}

\begin{abstract}  Let $\cM_{g,n}(\mathbf{L})$ be the moduli space of hyperbolic surfaces of genus $g$ with $n \geq 0$ hyperbolic ends of widths $\mathbf{L} \in \R_{\geq 0}^n$. We regard the total mass $|\mu_X^\kappa|$ of the Brownian loop measure with the killing rate $\kappa$ as a random variable on $\cM_{g,n}(\mathbf{L})$. Under the condition $|\mathbf{L}|^2 =o(g)$ as $g \to \infty$, we obtain the following two main results:
\ben
\item For any $\kappa > 0$, the expected value of $|\mu_X^\kappa|$ on all non-peripheral homotopy classes over $\cM_{g,n}(\mathbf{L})$ converges to an explicit function of $\kappa$, which blows up at the rate $ \log \left(\frac{1}{\kappa}\right)$ as $\kappa \to 0^+$.
\item For $\kappa=0$, over $\cM_{g,n}(\mathbf{L})$ the expected value of $|\mu_X|$ on homotopy classes of (iterates of) all non-peripheral simple closed geodesics is asymptotically $\frac{1}{2} \log g$.
\een
\end{abstract}

\tableofcontents

\section{Introduction}

    The Brownian loop measure was first introduced and studied on the complex plane $\C$ by Lawler, Schramm, and Werner; see e.g. \cite{Lawler_2003, lawler}. This measure admits a natural extension to arbitrary Riemann surfaces. Given a complete Riemann (hyperbolic) surface $(X,h)$, consider the Brownian motion generated by the Laplacian operator $\Delta_h$. The Brownian loop measure is then defined as
    \beq \mu_X := \int_0^\infty \frac{dt}{t} \int_X \W^t_{x \to x}(X) \dvol_h(x). \eeq
    This construction produces a measure on the set of oriented closed curves on $X$, modulo positive-oriented time parameterization. For Brownian motion with a constant killing rate $\kappa \geq -1/4$, the analogous formula defines a measure \(\mu_X^\kappa\), known as the Brownian loop measure with killing rate $\kappa$. For applications of the Brownian loop measure in the study of random conformal geometry, one may refer to  \cite{Wer08, lejan2011, AngPark, Wang_2025} and the references therein for related developments.

    Fix $n \geq 0$ and let $\mathbf{L} = (L_1,\cdots,L_n)\in \R_{\geq 0}^n$. We denote by $\cM_{g,n}(\mathbf{L})$ the moduli space of Riemann (hyperbolic) surfaces of genus $g$ with $n$ geodesic boundary components of lengths $L_1,\cdots,L_n$. When \(L_i = 0\) for some \(i \in \{1, \dots, n\}\), the corresponding boundary component is a cusp. 
    For any $X \in \cM_{g,n}(\mathbf{L})$, the surface $\bar X$ obtained by attaching funnels to all relevant boundary components of $X$ is a complete hyperbolic surface. The number of funnels of $\bar X$ equals that of the boundary components of $X$ with positive length.
    
    \begin{figure}[h]
        \centering
        \includegraphics[width=\linewidth]{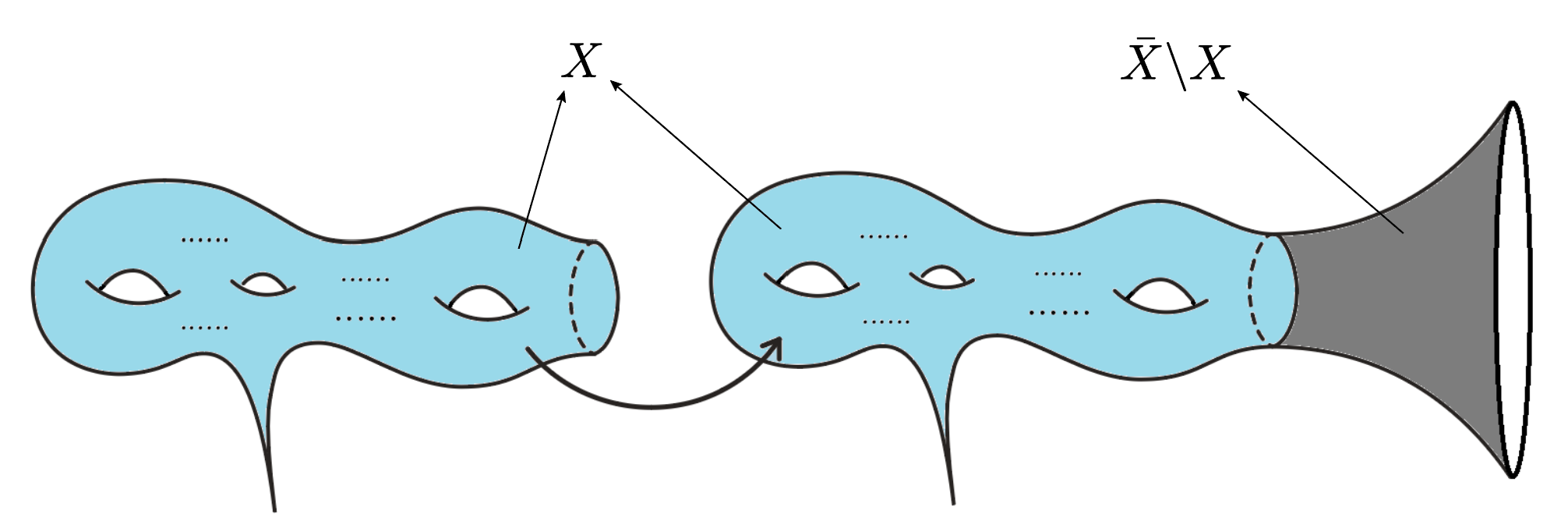}
        \caption{Attaching funnels to $X$. $\text{LHS} :X;\ \text{RHS} : \bar X$.}
        \label{attaching funnels}
    \end{figure}
    
    Choose $X\in  \cM_{g,n}(\mathbf{L})$. It is known from \cite{Wang_2025} of Wang-Xue and \cite{lemonde2026brownianloopsselbergzeta} of Lemonde-Wang that the Brownian loop measure $\mu_{\bar X}^\kappa$ of the set $\mathcal{C}_\alpha$ of all oriented closed curves homotopic to a given loop $\alpha \subset \overline{X}$ depends only on the length of the corresponding unique closed geodesic in $X$. In particular, if $\alpha$ is a closed geodesic on the boundary of $X$, then $\mu_{\bar X}^\kappa(\mathcal{C}_\alpha)$ depends solely on the length of $\alpha$, and not on the complex structure of $X$. 
    
    In this work, we study the Brownian loop measure of non-peripheral homotopy classes, i.e.  the ones that are neither bounding a disk or a cusp, nor iterates of a boundary geodesic component of $X$. We denote by $|\mu_{\bar X}^\kappa|$ the total mass of $\mu_{\bar X}^\kappa$ over all non-peripheral homotopy classes on $\bar X$. This in particular gives a random variable on $\cM_{g,n}(\mathbf{L})$.

    Denote by $\EgnL \Big[ |\mu_{\bar X}^\kappa| \Big]$ the expectation of $|\mu_{\bar X}^\kappa|$ over $\cM_{g,n}(\mathbf{L})$:
    \beq \EgnL \Big[ |\mu_{\bar X}^\kappa| \Big]:= \frac{1}{V_{g,n}(\mathbf{L})} \int_{\cM_{g,n}(\mathbf{L})} |\mu_{\bar X}^\kappa| dX, \eeq
    where $dX$ is the Weil-Petersson measure on $\cM_{g,n}(\mathbf{L})$ and $V_{g,n}(\mathbf{L}) < +\infty$ is the Weil-Petersson volume of $\cM_{g,n}(\mathbf{L})$.

    In this work, we investigate the asymptotic behavior of $\EgnL \Big[ |\mu_{\bar X}^\kappa| \Big]$ as $g \to \infty$. Our first result is as follows.
 \bt \label{Main theorem} Fix $\kappa > 0$ and $n \geq 0$. Let $\mathbf{L} =\mathbf{L}(g) \in \R_{\geq 0}^n$ satisfy $|\mathbf{L}|^2 = o(g)$ as $g \to \infty$. Then 
    \be \lim_{g \to \infty} \EgnL \Big[|\mu_{\bar X}^\kappa| \Big] = c \Bigg(\frac{1}{2} - \sqrt{\frac{1}{4} + \kappa} \Bigg), \ene
    where
    \beqar
c(z) &=& \log \frac{(z-1)^2}{z(z-2)} + 2S(z) - S(z+1) - S(z-1),\\
S(z) &=& \sum_{j=2}^\infty \frac{z^j}{j} \Big[\zeta(j+1) \zeta(j) - 1\Big],
    \eeqar
    and $\zeta$ denotes the Riemann zeta function. 
    \et
As $\kappa \to 0^+$, the limit $\lim\limits_{g \to \infty} \EgnL \Big[|\mu_{\bar X}^\kappa| \Big] \sim \log\left(\frac{1}{\kappa} \right)$.

     \begin{rem*} Let $X$ be a complete hyperbolic surface of finite area. Here, the case \(n=0\) yields a closed hyperbolic surface, whereas \(\mathbf{L}=0\) with \(n>0\) corresponds to a cusped hyperbolic surface. The Selberg Zeta function of $X$ is defined by
    \beq Z_X(s) := \prod_{\gamma \in \cP(X)} \prod_{k=0}^\infty \Big(1 - e^{(k + s)\ell_\gamma(X)} \Big), \eeq
    where $\cP(X)$ is the set of all primitive oriented closed geodesics in $X$. The function
    $Z_X(s)$ has meromorphic extension to $\C$, with a simple zero at $s = 1$ and no other zeros for $\mathrm{Re}(s) \geq 1$. It is known from \cite[Theorem 1]{lemonde2026brownianloopsselbergzeta} that for any $\kappa > 0$ and any such $X$,
    \beq |\mu_{X}^\kappa| = - \log Z_X \Bigg( \frac{1}{2} + \sqrt{\frac{1}{4} + \kappa} \Bigg). \eeq
    Thus, our main result Theorem \ref{Main theorem} for the finite-area case may be equivalently reformulated as follows: for any $\eps > 0$, 
    \be \lim_{g \to \infty} \mathbb{E}_{\mathrm{WP}}^{g,n} \Big[ - \log Z_X (1 + \eps) \Big] = c(-\eps). \ene
    \end{rem*}
    
    \vone
    
     For the Brownian loop measure without killing, i.e. $\kappa=0$, we investigate the asymptotic behavior of $\mu_{\bar X}$ over (iterates of) non-peripheral simple closed curves for large genus. Denote by $\cC^s(\bar X)$ the set of all loops on $\bar X$ such that each loop is freely homotopic to some iterate of a non-peripheral simple closed geodesic on $X$. Our second result is as follows.
    \bt \label{Total mass of simple curves, case kappa = 0} Fix $n \geq 0$ and let $\mathbf{L} = \mathbf{L}(g) \in \R_{\geq 0}^n$. If $|\mathbf{L}|^2 = o(g)$, then 
    \be \lim_{g \to \infty} \frac{\EgnL \Big[\mu_{\bar X}(\cC^s(\bar X)) \Big]}{\log g} = \frac{1}{2}.  \ene
    Moreover, if the stronger condition $|\mathbf{L}|^2 = o(\frac{g}{\log g})$ holds, then there is a constant $C_n$ that depends only on $n$ such that for any $g$ sufficiently large, 
    \be \Bigg| \EgnL \Big[\mu_{\bar X}(\cC^s(\bar X)) \Big]  - \frac{1}{2} \log g \Bigg| \leq C_n.  \ene
    \et

    \vone

    \subsection*{Notations.} 
    We write two positive functions $f_1(g)$ and $f_2(g)$  as $$f_1(g)\prec f_2(g) \ \ \text{ or } \quad f_2(g)\succ f_1(g)\ \ \text{ or } \quad f_1(g)=O(f_2(g))$$ if there exists a constant $C>0$ independent of $g$, such that $f_1(g) \leq C \cdot f_2(g)$. We say $f_1(g) \sim f_2(g)$ if $\lim\limits_{g \to \infty} \frac{f_1(g)}{f_2(g)}  = 1$, and  $f_1(g) =o( f_2(g))$ if 
    $\lim\limits_{g \to \infty} \frac{f_1(g)}{f_2(g)}  = 0$. For $\mathbf{L} = (L_1,\cdots,L_n) \in \R^n$, we denote the squared \(\ell^2\)-norm by $|\mathbf{L}|^2 = \sum_{i=1}^n L_i^2, $
    and its \(\ell^1\)-norm by
    $|\mathbf{L}|_1 = \sum_{i=1}^n |L_i|.$

\subsection*{Plan of the paper}
In Section \ref{sec-pre} we recall some basic properties of the Brownian loop measure on hyperbolic surfaces, two-dimensional hyperbolic geometry, and the Weil-Petersson geometry of moduli space of Riemann surfaces. Section \ref{sec3-scc} is devoted to the study of the Brownian loop measure of simple closed loops, culminating in the proof of Theorem \ref{Total mass of simple curves, case kappa = 0}. Finally, Section \ref{sec4-nscc} completes the proof of Theorem \ref{Main theorem} by additionally studying the Brownian loop measure of non-simple closed loops.

    \subsection*{Acknowledgements.} We are grateful to Yuxin He and Yuhao Xue for their valuable comments and suggestions on an earlier version of this manuscript. We also thank Roman Lemonde, Yilin Wang and Hao Wu for their interest in this work. The authors are partially supported by the National Key R \& D Program of China (2025YFA1017500) and NSFC grants No. 12361141813 and 12425107. 
    
    \vone

    \section{Preliminary}\label{sec-pre}
    
    In this section, we review the key notation and fundamental properties of the Brownian loop measure, hyperbolic surfaces, and the Weil-Petersson geometry of moduli spaces.

    \vone

    \subsection{The Brownian loop measure on hyperbolic surfaces} Let $(X, h)$ be an orientable Riemann surface equipped with a conformal metric $h$. Consider the Brownian motion on $X$ generated by the Beltrami-Laplace operator $\Delta_h$. The parametrized Brownian loop measure on $X$ is defined as a measure on the space of all loops $\{W \in C([0,t],X);  t \geq 0, W_0 = W_t\}$ given by
    \beq \mu_X^\ast = \int_0^\infty \frac{dt}{t} \int_X \W^t_{x \to x}(X) \dvol_h(x). \eeq
    The Brownian loop measure $\mu_X$ is the induced measure on the space of all loops, modulo orientation-preserving time parameterization. One remarkable property of $\mu_X$ is its conformal invariance (see e.g. \cite{lawler, Wer08, AngPark}).

    \bt Let $X_1 = (X,h)$ and $X_2 = (X,e^{2\sigma} h)$ be conformally equivalent Riemannian surfaces, then $\mu_{X_1} = \mu_{X_2}$.
    \et

    The uniformization theorem for Riemann surfaces states that any simply-connected Riemann surface is biholomorphic to one of $\S^2,\C$ or $\H$. Moreover, $\S^2$ covers only $\S^2$ and $\C$ covers only $\C,\C^\ast$ and $\T^2$. It follows that every Riemann surface $X$ other than these cases admits a complete hyperbolic metric compatible with its complex structure. For such surfaces, the Brownian loop measure $\mu_X$ can be defined using the complete hyperbolic metric on $X$ associated to the complex structure of $X$, allowing us to compute $\mu_X$ explicitly in terms of the hyperbolic length spectrum.

    Let $(X,h)$ be a complete hyperbolic surface, $\gamma$ be an oriented primitive closed geodesic on $X$, and $\cC_X(\gamma^m)$ be the free homotopy class of $\gamma^m$ on $X$, where $m \geq 1$. It is known from \cite[Lemma 3.2]{Wang_2025} of Wang-Xue that
    \be \mu_X(\cC_X(\gamma^m)) = \frac{1}{m} \frac{1}{e^{m \ell_\gamma(X)} - 1}, \ene
    where $\ell_\gamma(X)$ is the hyperbolic length of $\gamma$ on $X$. Moreover, if a closed curve $\gamma$ is either homotopically trivial or homotopic to a cusp, then $\mu_X(\cC_X(\gamma)) = +\infty$. A homotopy class $\cC$ is called \emph{non-peripheral} if it is none of the following: a trivial homotopy class, a class homotopic to a cusp, or a class homotopic to an iterate of the boundary geodesic of a funnel component of $X$ (if exists). Every non-peripheral homotopy class contains a unique closed geodesic. Let $\cP(X)$ denote the set of all non-peripheral oriented primitive closed geodesics on $X$. The total mass of $\mu_X$ over all non-peripheral homotopy classes is given by
    \be \sum_{m=1}^\infty  \sum_{\gamma \in \cP(X)} \mu_X(\cC_X(\gamma^m)) = \sum_{m=1}^\infty \sum_{\gamma \in \cP(X)}  \frac{1}{m} \frac{1}{e^{m \ell_\gamma(X)} - 1}. \ene
    We denote this summation by $|\mu_X|$. By \cite[Corollary 4.9]{Wang_2025} we know that $|\mu_X| < +\infty$ for any complete hyperbolic surface $X$ of finite type and infinite area.

    The definition of the Brownian loop measure extends naturally to  Brownian motions with a constant killing rate $\kappa \geq -\frac{1}{4}$, denoted by $\mu_X^\kappa$. If $(X,h)$ is a complete hyperbolic surface and $\gamma \in \cP(X),\ m \geq 1, \ \kappa \geq - \frac{1}{4}$, it is known from \cite[Lemma 3.1]{lemonde2026brownianloopsselbergzeta} of Lemonde-Wang that
    \be\label{LW26-0} \mu_X^\kappa(\cC_X(\gamma^m)) = \frac{1}{m} \frac{e^{k \cdot m \ell_\gamma(X)}}{e^{m \cdot \ell_\gamma(X)} - 1},\ \text{ where } k = \frac{1}{2} - \sqrt{\frac{1}{4} + \kappa}. \ene
    Therefore, the total mass of $\mu_X^\kappa$ over all non-peripheral homotopy classes is given by
    \be\label{LW26} \sum_{m=1}^\infty  \sum_{\gamma \in \cP(X)} \mu_X^\kappa(\cC_X(\gamma^m)) = \sum_{m=1}^\infty \sum_{\gamma \in \cP(X)} \frac{1}{m} \frac{e^{k \cdot m \ell_\gamma(X)}}{e^{m\ell_\gamma(X)} - 1}.\ene
    We denote this total mass by $|\mu_X^\kappa|$. If $X$ is geometrically finite and satisfies $\delta_X + k < 1$, where $\delta_X$ is the critical exponent of $X$, then \(|\mu_X^\kappa|\) is finite. For convenience, we state and prove the following fact:

    \bpro \label{Finiteness of total Brownian measure} Let $X$ be a geometrically finite, non-elementary complete hyperbolic surface, and $\delta$ be the critical exponent of $X$. Then the total mass of $\mu_X^\kappa$ over all non-peripheral homotopy classes is finite iff $\delta + k < 1.$
    \epro

    \bp The prime geodesic theorem (see e.g. \cite{Guillope1986, buser2010geometry, Borth}) for $X$ asserts that
    \beq N_X(L) =  \# \{\gamma \in \cP(X);\ \ell_\gamma(X) \leq L\} \sim \frac{e^{\delta L}}{\delta L},\ \text{ as } L \to \infty. \eeq
    We remark here that this asymptotic behavior remains valid even though our set \(\mathcal{P}(X)\) excludes peripheral simple closed geodesics of $X$. From \eqref{LW26}, the total mass of $\mu_X^\kappa$ over all non-peripheral homotopy classes is 
    \beq |\mu_X^\kappa| = \sum_{\gamma \in \cP(X)} \sum_{m=1}^\infty \frac{1}{m} \frac{e^{k \cdot m \ell_\gamma(X)}}{e^{m\ell_\gamma(X)} - 1}. \eeq
    Since $X$ contains only finitely many closed geodesics with length $< 1$, there is a constant $0<c_X<\infty$ such that
    \beq \sum_{\ell_\gamma(X) < 1} \sum_{m=1}^\infty \frac{1}{m} \frac{e^{k \cdot m \ell_\gamma(X)}}{e^{m\ell_\gamma(X)} - 1} = c_X < \infty. \eeq
    For any $\ell_\gamma(X) \geq 1$ and $m \geq 1$, we have $m(e^{m \ell_\gamma(X)} - 1) \geq \frac{e^{m \ell_\gamma(X)}}{2}$. Thus, there exists a constant $c_X>0$ depending only on $X$ such that
    \beq |\mu_X^\kappa| \leq c_X + 2\sum_{\ell_\gamma(X) \geq 1} \sum_{m=1}^\infty \frac{e^{k \cdot m \ell_\gamma(X)}}{e^{m \ell_\gamma(X)}} = c_X  + 2 \sum_{\ell_\gamma(X) \geq 1} \frac{e^{(k-1) \ell_\gamma(X)}}{1 - e^{(k-1) \ell_\gamma(X)}} . \eeq
    Since $k \leq \frac{1}{2}$, for any $\ell_\gamma(X) \geq 1$ we have $e^{(k-1)\ell_\gamma(X)} < e^{-1/2} < \frac{2}{3}$. It follows that 
    \beq |\mu_X^\kappa| \leq c_X + 6 \sum_{\ell_\gamma(X) \geq 1} e^{(k-1) \ell_\gamma(X)} \leq c_X + 6 \sum_{L = 1}^\infty N_X(L+1) \cdot e^{(k-1)L}. \eeq
    By the prime geodesic theorem, there exists a constant $c'_X > 0$ such that $N_X(L) \leq c_X' e^{\delta L}$, for all $L \geq 1$. If $\delta + k < 1$, we obtain
    \beq |\mu_X^\kappa| \leq c_X + 6 c_X' \sum_{L=1}^\infty e^{(\delta + k - 1)L} < +\infty. \eeq
    Conversely, $|\mu_X^\kappa|$ satisfies
    \beq |\mu_X^\kappa|>\sum_{\gamma \in \cP(X)} e^{(k-1) \ell_\gamma(X)} \geq  \sum\limits_{m=0}^\infty  e^{(k-1)m}\cdot \left(N_X(m+1)-N_X(m)\right).
\eeq
    By the prime geodesic theorem for $X$, given any $\eps > 0$, there is a constant $M>0$ such that for any $m\geq M$, $$ (1-\eps) \frac{e^{\delta m}}{\delta m} \leq N_X(m)  \leq (1+\eps) \frac{e^{\delta m}}{\delta m}.$$ 
    Consequently,
    \beq \begin{split} N_X(m+1) - N_X(m) & \geq (1-\eps) \frac{e^{\delta(m+1)}}{m+1} - (1+\eps) \frac{e^{\delta m}}{\delta m} \\ & = \Big( \frac{(1-\eps)m e^{\delta}}{m+1} - (1 + \eps) \Big) \cdot \frac{e^{\delta m}}{\delta m}. \end{split} \eeq
    If $\eps>0$ is sufficiently small and $M$ is sufficiently large, this yields
    \beq N_X(m+1) - N_X(m) \geq \frac{e^{\delta} - 1}{2} \cdot \frac{e^{\delta m}}{\delta m} \eeq
    for any $m \geq M$. If $\delta + k \geq 1$, it follows that
    \beq |\mu_X^\kappa|>\sum\limits_{m=M}^\infty  e^{(k-1)m}\cdot \frac{e^\delta - 1}{2} \frac{e^{\delta m}}{\delta m} = \infty. \eeq
    This completes the proof.
    \ep

    \bcor Let $X$ be a complete hyperbolic surface with $0<  |\chi(X)| < +\infty$, then the value of $|\mu_X^\kappa|$ has the following possibilities: \ben
        \item If $\area(X) < +\infty$ and $\kappa > 0$, then $|\mu_X^\kappa| < +\infty$.
        \item If $\area(X) < +\infty$ and $\kappa = 0$, then $|\mu_X^\kappa| = |\mu_X| = +\infty$.
        \item If $\area(X) = +\infty$ and $\kappa \geq 0$, then $|\mu_X^\kappa| < +\infty$.
    \een
    \ecor

    \bp By assumption, $X$ is geometrically finite and non-elementary. For such surfaces, $0 < \delta_X \leq 1$, with the equality $\delta_X = 1$ iff $X$ has finite area, see, e.g. \cite[Section 2.2]{Borth}. By definition, $\kappa\geq 0$ iff $k\leq 0$. The conclusion then follows directly from Lemma \ref{Finiteness of total Brownian measure}.
    \ep

    \vone

    \subsection{Counting geodesics} \label{filling construction}
    In the following, we recall some counting results  for closed geodesics on hyperbolic surfaces. A general upper bound is given by the following lemma: see \cite[Lemma 6.6.4]{buser2010geometry} and \cite[Lemma 3.7]{he2026uniformspectralgapsrandom} for more details.

    \bl \label{Classical bound of geodesics} Let $X \in \cM_{g,n}(\mathbf{L})$ and $L > 0$. Then there are at most $(g-1+\frac{n}{2}) \cdot e^{L + 6}$ oriented closed geodesics on $X$ of length $\leq L$, which are not iterates of simple closed geodesics with length $\leq 2\sinh^{-1}(1)$.
    \el

    \vone

    A closed geodesic $\gamma \subset X$ is called \emph{filling} if each component of $X \backslash \gamma$ retracts onto one of the following: a point, a cusp or a boundary geodesic component of $X$. The following counting for filling closed geodesics is established in \cite[Theorem 4]{Wu_2022}, also \cite[Theorem 18]{Wu_2025}:
    \bl \label{main counting WX22} For any $0 < \eps < \frac{1}{2}$ and $m = 2g-2+n \geq 1$, there exists a constant $C(\eps,m)$ such that for any hyperbolic surface $X \in T_g(x_1,\cdots,x_n)$,
    \beq \#_f(X,L) \leq C(\eps,m) \cdot e^{L - \frac{1-\eps}{2} \sum x_i}, \eeq
    where $\#_f(X,L)$ is the number of filling closed geodesics in $X$ of length $\leq L$.
    \el

    \vone

    Now, given $X \in \cM_{g,n}(\mathbf{L})$ and a non-simple closed geodesic $\gamma \subset X$, we construct a geodesic subsurface $X(\gamma)$ of $X$ as in e.g. \cite{Mirzakhani_2019, Wu_2022, Nie_2023}: Let $\cN_\eps(\gamma)$ be a small tubular neighborhood of $\gamma$. For any connected component $\xi$ of $\p \cN_\eps(\gamma)$, we deform $\cN_{\eps}(\gamma)$ into $X(\gamma)$ as follows: If $\xi$ is homotopically trivial, we fill the disk bounded by $\xi$ into $\cN_\eps(\gamma)$; if $\xi$ is homotopic to a cusp, then we fill the cusp bounded by $\xi$ into $\cN_\eps(\gamma)$. If none of these happen, then $\xi$ is homotopic to a unique simple closed geodesic, and we deform $\xi$ to this simple closed geodesic.

    \begin{figure}[h]
        \centering
        \includegraphics[width=0.8\linewidth]{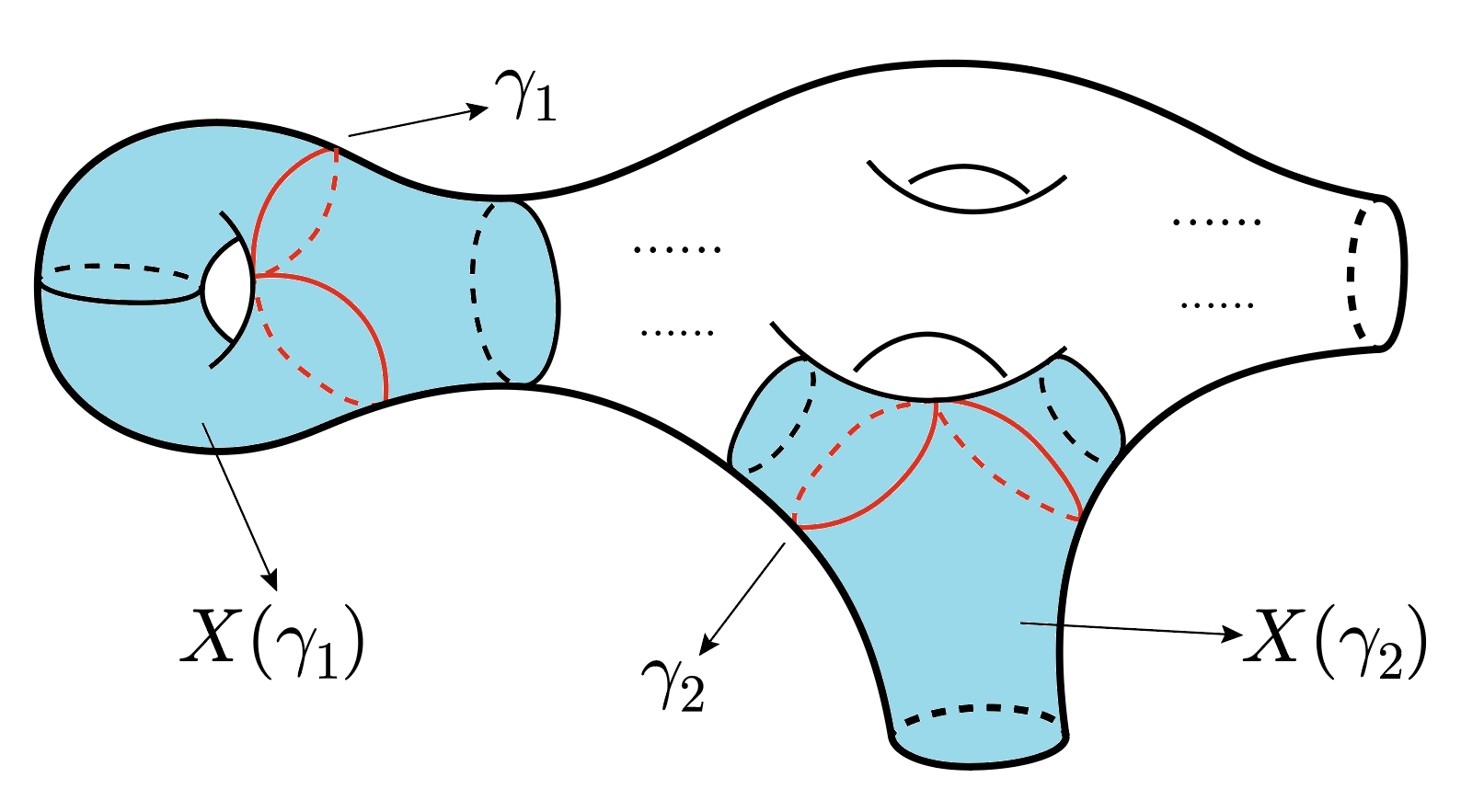}
        \caption{Examples of the filling construction}
        \label{fig:placeholder}
    \end{figure}
    
    We note, as in \cite{Wu_2022}, that if two components of $\cN_\eps(\gamma)$ deform to the same simple closed geodesic, they are not glued together, i.e. $X(\gamma)$ is regarded as an open subsurface of $X$.

    By construction, $\gamma$ is filling in $X(\gamma)$, and $\ell(\p X(\gamma)) \leq 2 \ell(\gamma)$. The standard isoperimetric inequality then implies that $\area(X(\gamma)) \leq 4 \ell(\gamma)$ (see also \cite[Proposition 7]{Wu_2022}). For $T > 0$, let $\mathrm{Sub}_T(X)$ denote the set of all subsurfaces $Y \subseteq X$ with geodesic boundary, satisfying $\area(Y) \leq 4T$ and $\ell_X(\p Y) \leq 2T$. In particular,  if $\ell_\gamma(X) \leq T$, then its associated filling subsurface $X(\gamma) \in \mathrm{Sub}_T(X)$.

    \vone

    \subsection{The Weil-Petersson metric} For $2g+n-2 \geq 1$ with $g,n \geq 0$ and $\mathbf{L} = (L_1,\cdots,L_n) \in \R_{\geq 0}^n$, let $\cT_{g,n}(\mathbf{L})$ denote the Teichm\"uller space of bordered hyperbolic surfaces of genus $g$ with $n$ geodesic boundary components of lengths $L_1,\cdots,L_n$. The mapping class group $\Mod_{g,n}$  acts on $\cT_{g,n}$ preserving the Weil-Petersson symplectic form $\omega_{\mathrm{WP}}$. The quotient space $\cM_{g,n}(\mathbf{L}) = \cT_{g,n}(\mathbf{L}) / \Mod_{g,n}$ is the moduli space of such Riemann surfaces. The Weil-Petersson volume form on $\cT_{g,n}(\mathbf{L})$ is given by
    \beq \dvol_{\mathrm{WP}} = \frac{1}{(3g-3+n)!} \omega_{\mathrm{WP}}^{3g-3+n}. \eeq
    This is a measure on $\cT_{g,n}(\mathbf{L})$ and invariant under $\mathrm{Mod}_{g,n}$, which descends to a measure on $\cM_{g,n}(\mathbf{L})$, denoted as $dX$ for brevity. Let $V_{g,n}(\mathbf{L})$ denote the total Weil-Petersson volume of $\cM_{g,n}(\mathbf{L})$. By \cite[Theorem 1.1]{Mirzakhani:2007a}, the volume $V_{g,n}(\mathbf{L})$ is a polynomial in $L_1,\cdots,L_n$ of degree $6g+2n-6$. We write $V_{g,n} = V_{g,n}(0)$. The Weil-Petersson probability measure on $\cM_{g,n}(\mathbf{L})$ is then defined by
    \beq \P(A) = \frac{1}{V_{g,n}(\mathbf{L})} \int_{\cM_{g,n}(\mathbf{L})} \mathbf{1}_A dX. \eeq

    The following integration formula of Mirzakhani is very useful for our calculations, one may see \cite[Theorem 7.1]{Mirzakhani:2007a} or \cite[Theorem 2.2]{Mirzakhani_2019} for more details.

    \bt[Mirzakhani] \label{Mirzakhani's integration formula} For any multi-curve $\Gamma = \cup_{i=1}^k \gamma_i$, where $\gamma_1,\cdots,\gamma_k$ are mutually disjoint, non-peripheral simple closed geodesics, there exists a constant $C_\Gamma \in (0,1]$ such that for any $F:\R^k \to \R_+$,
    \beq \begin{split} & \int_{\cM_{g,n}(\mathbf{L})} \sum_{(\alpha_1,\cdots,\alpha_k) \in \cO_\Gamma} F(\ell_{\alpha_1}(X),\cdots,\ell_{\alpha_k}(X)) dX 
    \\ & = C_\Gamma \int_{\R_+^k} F(x) V_{g,n}(\Gamma,x) \cdot  x_1 \cdots x_k dx_1 \cdots dx_k. \end{split} \eeq
    Here $\cO_\Gamma = \{(h \cdot \gamma_1,\cdots,h \cdot \gamma_k);\ h \in \Mod_{g,n}\}$ is the orbit of $\Gamma$ under the action of $\Mod_{g,n}$, and $V_{g,n}(\Gamma,x)$ is the Weil-Petersson volume of the moduli space of Riemann surfaces homeomorphic to $S_{g,n} \backslash \Gamma$ with boundary lengths given by $x,\mathbf{L}$. Moreover, $C_\Gamma = \frac{1}{2}$ if $g >2$ and $\Gamma$ is a simple non-separating closed curve.
    \et

    \vone

    \subsection{Bounds on Weil-Petersson volumes} In this subsection, we recall some well-known bounds on Weil-Petersson volumes, which will be used repeatedly in the proofs of Theorems \ref{Main theorem} and \ref{Total mass of simple curves, case kappa = 0}. It is known from \cite[Lemma 3.2]{Mir13a} that
    \be V_{g,n} \leq V_{g,n}(x_1,\cdots,x_n) \leq e^{(x_1 + \cdots + x_n)/2} V_{g,n}, \label{Mirzakhani's volume comparison formula}  \ene
    and that for any fixed $n \geq 0$, as $g \to \infty$
    \be \frac{V_{g,n}}{V_{g-1,n+2}} = 1 + O \Big( \frac{1}{g} \Big). \label{Volume asymptotics for adjacent g,n} \ene
    It is known from e.g. \cite[Proposition 3.1]{Mirzakhani_2019}, \cite[Lemma 22]{Nie_2023} or \cite[Lemma B.1]{he2026uniformspectralgapsrandom} that
    \be \Big(1 - c \cdot n \frac{\sum x_i^2}{g} \Big) \cdot \prod_{i=1}^n \frac{\sinh (x_i/2)}{x_i/2}  \leq \frac{V_{g,n} (x_1,\cdots,x_n)}{V_{g,n}} \leq \prod_{i=1}^n \frac{\sinh (x_i/2)}{x_i/2}, \label{NWX volume comparison} \ene
    where $c>0$ is a universal constant. Moreover, by \cite[Theorem 6]{Wu_2025}
    \be \frac{V_{g,n}(x_1,\cdots,x_n)}{V_{g,n}} \leq \Big(\prod_{i=1}^n \min \Big\{ \alpha(k) \cdot \frac{2g+n-2}{x_i^2}, 1 \Big\} \Big)^k \cdot \prod_{i=1}^n \frac{\sinh(x_i/2)}{x_i/2}, \label{New volume estimate} \ene
    for any $k \geq 1$, where $\alpha(k) > 0$ is a constant depending only on $k$. 
   
    For $r \geq 1$, following the notation in \cite{Nie_2023}, define 
    \beq W_r \overset{\text{def}}{=} \begin{cases}
        V_{\frac{r}{2} + 1} &\text{ if $r$ is even}, \\ \\ V_{\frac{r+1}{2},1} &\ \text{if $r$ is odd}.
    \end{cases} \eeq
    By \cite[Lemma 24]{Nie_2023}, there exists two universal constants $c,D>0$, such that for any $q \geq 1, n_1,\cdots,n_q \geq 0$ and $r \geq 2$, 
    \be \sum_{\{g_i\}} V_{g_1,n_1} \cdots V_{g_q,n_q} \leq c \Big(\frac{D}{r}\Big)^{q-1} W_r, \label{NWX volume summation estimate} \ene
    where the summation is taken over all tuples $\{g_i\}_{i=1}^q \subset \N^q$ satisfying 
    \beq \sum_{i=1}^q (2g_i-2+n_i) = r, \ 2g_i - 2 + n_i \geq 1. \eeq
    By \cite[Proposition 25]{Nie_2023} (or see \cite[Lemma B.17]{he2026uniformspectralgapsrandom} for a more general form), given $n \geq 0$ and $m \geq 1$, there exists a constant $c(m,n) > 0$, such that for any $g \geq m+1$, $q \geq 1$, and $n_1,\cdots,n_q \geq 1$, 
    \be \sum_{\{g_i\}} V_{g_1,n_1} \cdots V_{g_q,n_q} \leq c(m,n) \frac{V_{g,n}}{g^m}, \label{Euler char lowers by m} \ene
    where the summation is taken over all tuples $\{g_i\}_{i=1}^q \subset \N^q$ satisfying
    \beq \sum_{i=1}^q (2g_i-2+n_i) = 2g + n -2-m,\ 2g_i - 2 + n_i \geq 1. \eeq
    We conclude this section with two simple corollaries of these volume estimates, which will be used repeatedly in Section  \ref{sec3-scc}.

    \bl \label{Volume bound I} Let $n \geq 0$. There exists a constant $c_n>0$, depending only on $n$, such that for any $\mathbf{L}$ with $|\mathbf{L}|^2 = o(g)$ as $g \to \infty$ and all sufficiently large $g$, 
    \beq \prod_{i=1}^n \frac{\sinh(L_i/2)}{L_i/2} \cdot \sum_{g_1,n_1} V_{g_1,n_1+1} V_{g-g_1,n+1-n_1} \leq c_n \frac{V_{g,n}(\mathbf{L})}{g}, \eeq
    where the summation runs over all $g_1, n_1\geq 0$ satisfying $g_1 \leq g, \ n_1 \leq n$, and
    \beq 2g_1 - 1 + n_1 \geq 1,\ 2(g - g_1) + (n - n_1 - 1) \geq 1. \eeq
    \el

    \bp By \eqref{NWX volume summation estimate}, there exists a universal constant $c>0$ such that 
    \beq \sum_{g_1,n_1} V_{g_1,n_1+1} V_{g-g_1,n+1-n_1} \leq \sum_{n_1=0}^{n} \frac{c \cdot W_{2g+n-2}}{2g+n-2}. \eeq
    As $g \to \infty$, relation \eqref{Volume asymptotics for adjacent g,n} gives $W_{2g+n-2} \prec V_{g,n}$. Thus, there is a constant $c_n>0$ depending only on $n$ such that for sufficiently large $g$,
    \be  \sum_{g_1,n_1} V_{g_1,n_1+1} V_{g-g_1,n+1-n_1} \leq c_n \cdot \frac{V_{g,n}}{g}. \ene
    By \eqref{NWX volume comparison}, whenever $|\mathbf{L}|^2 = o(g)$, we have
    \beq \lim_{g \to \infty} \prod_{i=1}^n \frac{\sinh(L_i/2)}{L_i/2} \cdot \frac{V_{g,n}}{V_{g,n}(\mathbf{L})} = 1. \eeq
    Therefore, replacing \(c_n\) by \(2c_n\), we obtain for all sufficiently large $g$,
    \beq  \prod_{i=1}^n \frac{\sinh(L_i/2)}{L_i/2} \cdot \sum_{g_1,n_1} V_{g_1,n_1+1} V_{g-g_1,n+1-n_1} \leq c_n \frac{V_{g,n}(\mathbf{L})}{g}. \eeq
    This completes the proof.
    \ep

    \bl \label{Volume bound II} Let $n \geq 0$. There exists a constant $c_n>0$, depending only on $n$, and a universal constant $\alpha > 0$, such that the following holds: \ben
    \item For any $g \geq 2c_n(1 + |\mathbf{L}|^2)$ and any $0 <x < \sqrt{\frac{g}{4c_n}}$,
    \beq \begin{split} \Big| \frac{V_{g-1,n+2}(\mathbf{L},x,x)}{V_{g,n}(\mathbf{L})} &- \Big( \frac{\sinh(x/2)}{x/2} \Big)^2 \Big| \\ & \leq \Big( \frac{\sinh(x/2)}{x/2} \Big)^2 \cdot \frac{3c_n(1 + |\mathbf{L}|^2 + x^2)}{g}. \end{split} \eeq
    \item For any $g\geq 2c_n(1 + |\mathbf{L}|^2)$ and any $x > 0$, 
    \beq \frac{V_{g-1,n+2}(\mathbf{L},x,x)}{V_{g,n}(\mathbf{L})}  \leq \min \Big\{ \Big( \alpha \cdot \frac{2g+n-2}{x^2} \Big)^2,\ 4 \Big\} \cdot  \Big( \frac{\sinh(x/2)}{x/2} \Big)^2. \eeq
    \een
    In particular, if $|\mathbf{L}|^2 = o(g)$, then for any fixed $x > 0$,
    \beq  \lim_{g \to \infty} \frac{V_{g-1,n+2}(\mathbf{L},x,x)}{V_{g,n}(\mathbf{L})} = \Big( \frac{\sinh(x/2)}{x/2} \Big)^2. \eeq
    \el

    \bp Using \eqref{Volume asymptotics for adjacent g,n}, \eqref{NWX volume comparison} and \eqref{New volume estimate}, there is a constant $c_n>0$ depending only on $n$ such that for all relevant $g,x,\mathbf{L}$, the following bounds hold:
    \beq 1 - \frac{c_n}{g} \leq \frac{V_{g-1,n+2}}{V_{g,n}} \leq 1 + \frac{c_n}{g}, \eeq
    
    \beq \Big( \frac{\sinh(x/2)}{x/2} \Big)^2 \prod_{i=1}^n \frac{\sinh(L_i/2)}{L_i/2} \cdot \Big(1 - c_n \frac{|\mathbf{L}|^2 + 2x^2}{g} \Big) \leq \frac{V_{g-1,n+2} (\mathbf{L},x,x)}{V_{g-1,n+2}}, \eeq
    
    \beq \frac{V_{g-1,n+2} (\mathbf{L},x,x)}{V_{g-1,n+2}} \leq \Big( \frac{\sinh(x/2)}{x/2} \Big)^2 \prod_{i=1}^n \frac{\sinh(L_i/2)}{L_i/2},  \eeq
    
    \beq \prod_{i=1}^n \frac{\sinh(L_i/2)}{L_i/2} \cdot \Big(1 - c_n \frac{|\mathbf{L}|^2}{g} \Big) \leq \frac{V_{g,n} (\mathbf{L})}{V_{g,n}} \leq \prod_{i=1}^n \frac{\sinh(L_i/2)}{L_i/2}, \eeq
    
    \beq \frac{V_{g-1,n+2} (\mathbf{L},x,x)}{V_{g-1,n+2}} \leq \Big( \alpha \cdot \frac{2g+n-2}{x^2} \Big)^2 \cdot \Big( \frac{\sinh(x/2)}{x/2} \Big)^2 \prod_{i=1}^n \frac{\sinh(L_i/2)}{L_i/2}. \eeq
    More precisely, the first bound follows from \eqref{Volume asymptotics for adjacent g,n}, the last one follows from \eqref{New volume estimate} and the remaining bounds follow from \eqref{NWX volume comparison}. We now estimate the following volume quotient by factoring it as:
    \beq \frac{V_{g-1,n+2}(\mathbf{L},x,x)}{V_{g,n}(\mathbf{L})} = \frac{V_{g-1,n+2}(\mathbf{L},x,x)}{V_{g-1,n+2}} \frac{V_{g-1,n+2}}{V_{g,n}} \frac{V_{g,n}}{V_{g,n}(\mathbf{L})}, \eeq
    If $g > 2 c_n (1 + |\mathbf{L}|^2)$, then $1 - c_n \frac{|\mathbf{L}|^2}{g} \geq \frac{1}{2}$. Combining the first and the last three of the five bounds above, we obtain for any $g > 2 c_n(1 + |\mathbf{L}|^2)$ and $x > 0$,
    \beq \frac{V_{g-1,n+2}(\mathbf{L},x,x)}{V_{g,n}(\mathbf{L})} \leq \min \Big\{ \Big( 2\alpha \cdot \frac{2g+n-2}{x^2} \Big)^2,\ 4 \Big\} \cdot \Big( \frac{\sinh(x/2)}{x/2} \Big)^2, \eeq
    which completes the proof of part (2). 

    Moreover, for all $g > 2 c_n(1 + |\mathbf{L}|^2)$ and $x < \sqrt{\frac{g}{4c_n}}$, we have 
    \beq 1 - c_n \frac{|\mathbf{L}|^2 + 2x^2}{g} > 0 \ \ \text{ and } \ \ 1 - \frac{c_n}{g} > 0. \eeq
    Combining the first four of the five bounds listed above, we get
    \beq \frac{V_{g-1,n+2}(\mathbf{L},x,x)}{V_{g,n}(\mathbf{L})} \leq \Big( \frac{\sinh(x/2)}{x/2} \Big)^2 \cdot\Big(1 + \frac{2c_n(1 + |\mathbf{L}|^2)}{g} \Big), \eeq
    \beq \frac{V_{g-1,n+2}(\mathbf{L},x,x)}{V_{g,n}(\mathbf{L})} \geq \Big( \frac{\sinh(x/2)}{x/2} \Big)^2 \cdot \Big(1 - \frac{c_n(1 + |\mathbf{L}|^2 + 2x^2)}{g} \Big). \eeq
    Therefore, for all $g \geq 2c_n(1 + |\mathbf{L}|^2)$ and any $x < \sqrt{\frac{g}{4c_n}}$,
    \beq \begin{split} \Big| \frac{V_{g-1,n+2}(\mathbf{L},x,x)}{V_{g,n}(\mathbf{L})} &- \Big( \frac{\sinh(x/2)}{x/2} \Big)^2 \Big| \\ & \leq \Big( \frac{\sinh(x/2)}{x/2} \Big)^2 \cdot \frac{3c_n(1 + |\mathbf{L}|^2 + x^2)}{g}, \end{split} \eeq
    which completes the proof of part (1).
\ep

\vone

    \section{The total mass of $\mu_{\bar X}^\kappa$ over simple closed curves} \label{sec3-scc}

    Recall from the introduction that we denoted by $\cC^s(\bar X)$ the set of all loops on $\bar X$ such that each loop is freely homotopic to some iterate of a non-peripheral simple closed geodesic in $X$. The main results of this section are as follows:

    \bt \label{main theorem simple curves}
    Fix $\kappa > 0$ and $n \geq 0$. If $|\mathbf{L}|^2 = o(g)$, with the same function $c$ as in Theorem \ref{Main theorem}, then
    \beq \lim_{g \to \infty} \EgnL \Big[\mu_{\bar X}^\kappa(\cC^s(\bar X))\Big]  = c \Bigg(\frac{1}{2} - \sqrt{\frac{1}{4} + \kappa} \Bigg). \eeq
    \et

    \bt [= Theorem \ref{Total mass of simple curves, case kappa = 0}] \label{Total mass of simple curves, case kappa = 0, duplicate} Fix $n \geq 0$. If $|\mathbf{L}|^2 = o(g)$, we have
    \be \lim_{g \to \infty} \frac{\EgnL \Big[\mu_{\bar X}(\cC^s(\bar X)) \Big] }{\frac{1}{2}\log g} = 1.  \ene
    Moreover, if $|\mathbf{L}|^2 = o(\frac{g}{\log g})$, then there is a constant $c_n>0$ that depends only on $n$ such that for any $g$ sufficiently large, we have
    \be \Bigg| \EgnL \Big[\mu_{\bar X}(\cC^s(\bar X))\Big]   - \frac{1}{2} \log g \Bigg| \leq c_n.  \ene
    \et

    \vone

    \subsection{Estimates of $\mu_{\bar X}^\kappa(\cC^s_{sep}(\bar X))$ and $\mu_{\bar X}^\kappa(\cC^s_{nsep}(\bar X))$} We first set up some basic notation. Given $X \in \cM_{g,n}(\mathbf{L})$, let $\cP(X)$ denote the set of all non-peripheral primitive oriented closed geodesics in $X$. We define \ben
    \item $\cP^s_{sep}(X) := \{\gamma \in \cP(X),\text{ $\gamma$ is simple and separating} \}$.
    \item $\cP^s_{nsep}(X) := \{\gamma \in \cP(X), \text{ $\gamma$ is simple and non-separating}\}$.
    \item $\cP^{ns}(X) := \{\gamma \in \cP(X), \text{ $\gamma$ is non-simple} \}$.
    \een

    \begin{figure}[h]
        \centering \hspace*{-1cm}
        \includegraphics[width=0.8\linewidth]{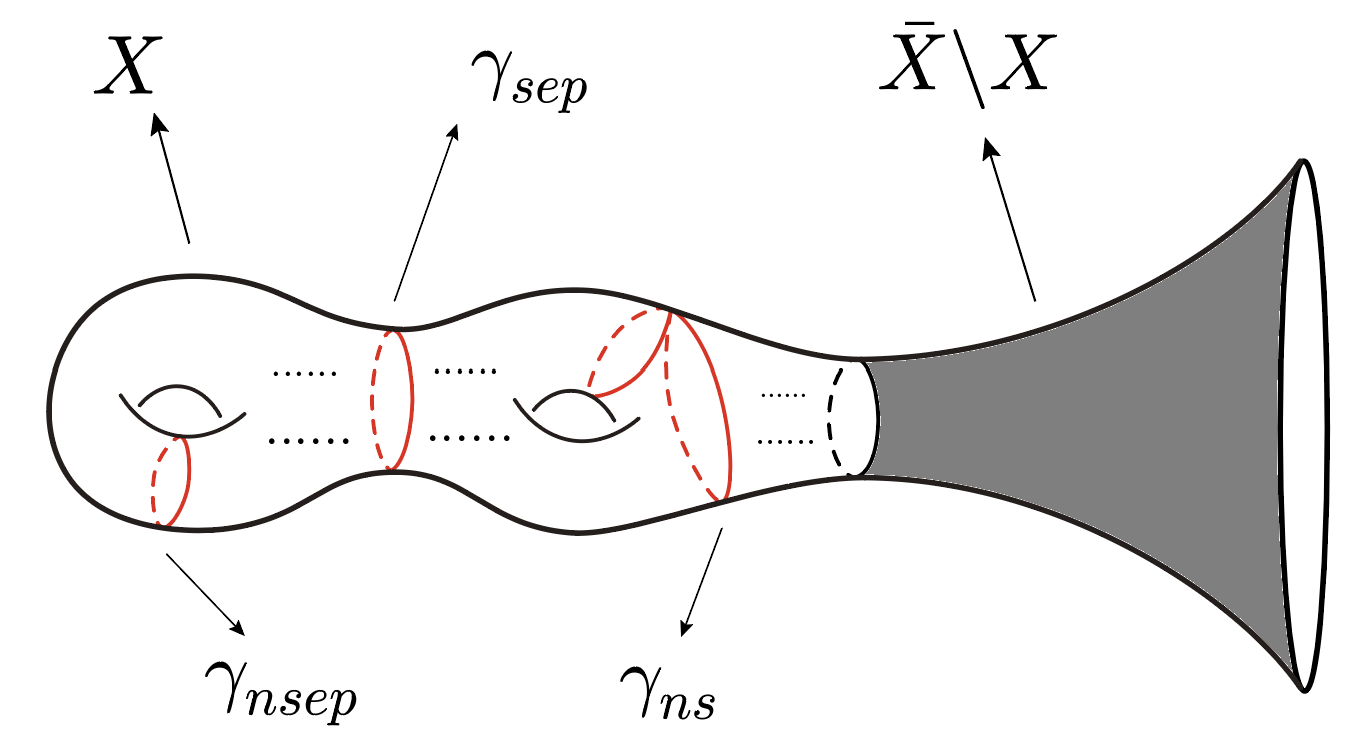}
        \caption{Examples: simple separating geodesics ($\gamma_{sep}$), simple non-separating geodesics ($\gamma_{nsep}$) and nonsimple geodesics ($\gamma_{ns}$) in $X$.}
        \label{fig:placeholder}
    \end{figure}
    
    \noindent Clearly, we have the disjoint union
    \beq \cP(X) = \cP^s_{sep}(X) \cup \cP^s_{nsep}(X) \cup \cP^{ns}(X). \eeq

    \noindent Recall that \(\bar X\) is obtained by attaching funnels to $X$ (when needed) to form a complete hyperbolic surface. We classify all the non-peripheral loops on $\bar X$ into three cases: \ben
    \item $\cC^s_{sep}(\bar X) = \{\text{loops homotopic to $\gamma^m$ for some $\gamma \in \cP^s_{sep}(X), m \geq 1$}\}$.
    \item $\cC^s_{nsep}(\bar X) = \{\text{loops homotopic to $\gamma^m$ for some $\gamma \in \cP^s_{nsep}(X), m \geq 1$}\}$.
    \item $\cC^{ns}(\bar X) = \{\text{loops homotopic to $\gamma^m$ for some $\gamma \in \cP^{ns}(X), m \geq 1$}\}$.
    \een
    By the definition of $|\mu_{\bar X}^\kappa|$, we have
    \beq \begin{split} |\mu_{\bar X}^\kappa| &= \mu_{\bar X}^\kappa(\cC^s_{sep}(\bar X)) + \mu_{\bar X}^\kappa(\cC^s_{nsep}(\bar X)) + \mu_{\bar X}^\kappa(\cC^{ns}(\bar X)) \\ & = \mu_{\bar X}^\kappa(\cC^s(\bar X)) + \mu_{\bar X}^\kappa(\cC^{ns}(\bar X)). \end{split}\eeq
    It is known from \cite[Lemma 3.2]{Wang_2025} of Wang-Xue and \cite[Lemma 3.1]{lemonde2026brownianloopsselbergzeta} of Lemonde-Wang that
    \be \label{Total mass of separating} \mu^\kappa_{\bar X}(\cC^s_{sep}(\bar X)) = \sum_{m = 1}^\infty \sum_{\gamma \in \cP^s_{sep}(X)} \frac{1}{m} \frac{e^{k \cdot m \ell_\gamma(X)}}{e^{m \ell_\gamma(X)} - 1}, \ene
    \be \label{Total mass of non-separating} \mu^\kappa_{\bar X}(\cC^s_{nsep}(\bar X)) = \sum_{m = 1}^\infty \sum_{\gamma \in \cP^s_{nsep}(X)} \frac{1}{m} \frac{e^{k \cdot m \ell_\gamma(X)}}{e^{m \ell_\gamma(X)} - 1}, \ene
    \be \label{Total mass of nonsimple} \mu^\kappa_{\bar X}(\cC^{ns}(\bar X)) = \sum_{m = 1}^\infty \sum_{\gamma \in \cP^{ns}(X)} \frac{1}{m} \frac{e^{k \cdot m \ell_\gamma(X)}}{e^{m \ell_\gamma(X)} - 1}. \ene
    Where $k = \frac{1}{2} - \sqrt{\frac{1}{4} + \kappa}$. These formulas for the mass of Brownian loop measure $\mu_{\bar X}^\kappa$ will be used in this and the next section to prove Theorem \ref{Main theorem} and Theorem \ref{Total mass of simple curves, case kappa = 0}.
    
    \bl \label{sep, m geq 2, k leq 0} Fix $n \geq 0$, and assume that $|\mathbf{L}|^2 = o(g)$ as $g \to \infty$. Then
    \beq \frac{1}{V_{g,n}(\mathbf{L})} \int_{\cM_{g,n}(\mathbf{L})} \sum_{m \geq 2} \sum_{\gamma \in \cP^s_{sep}(X)} \frac{1}{m(e^{m \ell_\gamma(X)} - 1)} dX \prec \frac{1}{g}. \eeq
    \el

    \bp Set the function $F$ as
    \beq F(x) = \sum_{m \geq 2} \frac{1}{m(e^{mx} - 1)}. \eeq
    For any $\gamma \in \cP^s_{sep}(X)$, the curve $\gamma$ separates $X$ into two subsurfaces with the following  topology:
    \be X \backslash \gamma \cong  S_{g_1,n_1+1} \bigcup S_{g - g_1,n + 1 - n_1},\ \textit{where } 0 \leq n_1 \leq \big[ \frac{n}{2} \big]. \label{topological type, case 1} \ene
    Let $\cS(g_1,n_1)$ denote the collection of all simple closed curves $\gamma$ satisfying the topological decomposition  \eqref{topological type, case 1}. Applying Mirzakhani's integration formula (Theorem \ref{Mirzakhani's integration formula}) to $F$, we obtain
    \beq \int_{\cM_{g,n}(\mathbf{L})} \Big[ \sum_{\gamma \in \cP^s_{sep}(X)} F(\ell_\gamma(X)) \Big] dX = \sum_{g_1,n_1} \sum_{[\gamma]} 2C_{\gamma}\int_0^\infty F(x) V_{g,n}(\gamma,x) \cdot xdx, \eeq
    where $[\gamma]$ runs over all equivalent classes in $\cS(g_1,n_1)/\mathrm{Mod}_{g,n}$. For each $[\gamma]$, the constant $C_\gamma \in (0,1]$ , the indices $g_1,n_1$ are determined via \eqref{topological type, case 1}, and 
    \beq V_{g,n}(\gamma,x) =  V_{g_1,n_1+1 }(\mathbf{L}_1,x) \cdot V_{g-g_1,n+1-n_1}(\mathbf{L}_2,x), \eeq where $\mathbf{L}_1 \sqcup \mathbf{L}_2 = \mathbf{L}$ denotes the partition of the boundary lengths of $X$, according to whether the corresponding boundary component lies on $S_{g_1,n_1 + 1}$ or $S_{g - g_1,n+1 - n_1}$. It is straightforward to see that for any $(g_1,n_1)$, there are at most $2^n$ choices for $[\gamma]$. For each $[\gamma]$, the volume bound \eqref{NWX volume comparison} yields
    \beq V_{g,n}(\gamma,x) \leq e^{x} \cdot \prod_{i=1}^n \frac{\sinh(L_i/2)}{L_i/2} \cdot V_{g_1,n_1+1} V_{g-g_1,n+1-n_1}. \eeq
    Substituting this bound into the integral, we get 
    \beq \begin{split} \int_{\cM_{g,n}(\mathbf{L})} \Big[ \sum_{\gamma \in \cP^s_{sep}(X)} F(\ell_\gamma(X)) \Big] dX & \leq 2^{n+1}  \prod_{i=1}^n \frac{\sinh(L_i/2)}{L_i/2} \int_0^\infty xF(x) e^xdx \\ & \times \sum_{g_1,n_1} V_{g_1,n_1+1} V_{g-g_1,n+1-n_1}. \end{split} \eeq
    Since $x F(x) e^x \in L^1(0,+\infty)$, the integral is bounded by a universal constant. Applying Lemma \ref{Volume bound I} to the remaining sum, we conclude
    \beq \int_{\cM_{g,n}(\mathbf{L})} \Big[ \sum_{\gamma \in \cP^s_{sep}(X)} F(\ell_\gamma(X)) \Big] dX \prec \frac{V_{g,n}(\mathbf{L})}{g}. \eeq
Dividing both sides by \(V_{g,n}(\mathbf{L})\) completes the proof.
    \ep

    \bl \label{sep, m = 1, k leq 0} Fix $n \geq 0$, and assume that $|\mathbf{L}|^2 = o(g)$ as $g \to \infty$. Then
    \beq \frac{1}{V_{g,n}(\mathbf{L})} \int_{\cM_{g,n}(\mathbf{L})} \sum_{\gamma \in \cP^s_{sep}(X)} \frac{1}{e^{ \ell_\gamma(X)} - 1} dX  \prec \frac{\log g}{g}. \eeq
    \el

    \bp Set the integration function $F$ as
    \beq F(x) = \frac{1}{e^x - 1}. \eeq
    Applying Mirzakhani's integration formula (Theorem \ref{Mirzakhani's integration formula}) to $F$, we obtain
    \beq \int_{\cM_{g,n}(\mathbf{L})} \Big[ \sum_{\gamma \in \cP^s_{sep}(X)} F(\ell_\gamma(X)) \Big] dX = \sum_{g_1,n_1} \sum_{[\gamma]} 2C_{\gamma}\int_0^\infty F(x) V_{g,n}(\gamma,x) \cdot xdx. \eeq
    By the volume estimate \eqref{New volume estimate}, there exists a universal constant $\alpha>0$ such that 
    \beq \begin{split} V_{g,n}(\gamma,x) &\leq \alpha \cdot \Big( \frac{\sinh(x/2)}{x/2} \Big)^2 \min\Bigg\{\Big(\frac{2g + n - 2}{x^2} \Big)^2, 1\Bigg\} \\ & \times \prod_{i=1}^n \frac{\sinh(L_i/2)}{L_i/2} \cdot V_{g_1,n_1+1} V_{g-g_1, n+1 - n_1}. \end{split} \eeq
    Substituting this bound into the integration formula and using Lemma \ref{Volume bound I}, we find a constant $c_n>0$ depending only on $n$ such that for any sufficiently large $g$, 
    \be \begin{split}
    \int_{\cM_{g,n}(\mathbf{L})} &\Big[ \sum_{\gamma \in \cP^s_{sep}(X)} F(\ell_\gamma(X)) \Big] dX  \leq c_n \frac{V_{g,n}(\mathbf{L})}{g} \\ & \times \int_0^\infty xF(x) \Big( \frac{\sinh(x/2)}{x/2} \Big)^2 \min \Big( \frac{g^2}{x^4}, 1 \Big) dx. \end{split} \ene
    Now we split the integral into two parts: 
    \beq \mathrm{I} =  \int_0^{\sqrt{g}} x F(x) \Big( \frac{\sinh(x/2)}{x/2} \Big)^2 \min \Big( \frac{g^2}{x^4}, 1 \Big) dx \eeq
    and
    \beq \mathrm{II} =  \int_{\sqrt{g}}^{\infty} x F(x) \Big( \frac{\sinh(x/2)}{x/2} \Big)^2 \min \Big( \frac{g^2}{x^4}, 1 \Big) dx. \eeq
    It follows that
    \beq \frac{1}{V_{g,n}(\mathbf{L})} \int_{\cM_{g,n}(\mathbf{L})} \sum_{\gamma \in \cP^s_{sep}(X)} \frac{1}{e^{ \ell_\gamma(X)} - 1} dX \leq c_n \frac{\mathrm{I} + \mathrm{II}}{g}. \eeq
    we estimate $\mathrm{I}$ and $\mathrm{II}$ separately. For $\mathrm{I}$, we have
    \beq \mathrm{I} \leq \int_{0}^{\sqrt{g}} x F(x) \Big( \frac{\sinh(x/2)}{x/2} \Big)^2 dx = \int_{0}^{\sqrt{g}} \frac{1 - e^{-x}}{x} dx \leq 1 + \frac{1}{2} \log g. \eeq
    For $\mathrm{II}$, we have
    \beq \mathrm{II} \leq \int_{\sqrt{g}}^\infty x F(x) \Big( \frac{\sinh(x/2)}{x/2} \Big)^2 \frac{g^2}{x^4}dx = g^2 \int_{\sqrt{g}}^\infty \frac{1 - e^{-x}}{x^5} dx \leq \frac{1}{4}. \eeq
    Combining the bounds for $\mathrm{I}$ and $\mathrm{II}$, we conclude
    \beq \frac{1}{V_{g,n}(\mathbf{L})} \int_{\cM_{g,n}(\mathbf{L})} \sum_{\gamma \in \cP^s_{sep}(X)} \frac{1}{e^{ \ell_\gamma(X)} - 1} dX  \prec \frac{\log g}{g}. \eeq
    This completes the proof.
    \ep

    \bpro \label{sep, m = all, k = 0} Fix $n \geq 0$, $\kappa \geq 0$, and assume that $|\mathbf{L}|^2 = o(g)$. Then there exists a constant $c_n>0$ depending only on $n$ such that for any sufficiently large $g$,
    \beq \EgnL \Big[\mu_{\bar X}^\kappa (\cC^{s}_{sep}(\bar X)) \Big] \leq c_n \frac{\log g}{g}. \eeq
    \epro

    \bp Recall that $\kappa \geq 0$ iff $k\leq 0$. The result then follows immediately from Lemmas \ref{sep, m geq 2, k leq 0} and \ref{sep, m = 1, k leq 0}.
    \ep
    
   Next, we treat the simple non-separating case.  
    \bpro \label{nsep, m=all, k<0} Fix $n \geq 0$, and assume $|\mathbf{L}|^2 = o(g) $. If $\kappa > 0$, then
    \beq \lim_{g \to \infty} \EgnL \Big[ \mu_{\bar X}^\kappa(\cC^{s}_{nsep}(\bar X)) \Big] = c(k), \eeq
    where $c(k)$ is a constant depending on $k = \frac{1}{2} - \sqrt{\frac{1}{4} + \kappa}$, given by
    \beq c(k) = \int_0^\infty \sum_{m=1}^\infty \frac{e^{k \cdot x}}{m(e^{mx} - 1)} \cdot \frac{(e^x - 1)^2}{x e^x} dx. \eeq
    \epro

    \bp Let $F(x)$ be the integration function defined by the series
    \beq F(x) = \sum_{m=1}^\infty \frac{e^{kx}}{m(e^{mx}-1)}. \eeq
    By identity \eqref{Total mass of non-separating}, we have
    \beq \mu_{\bar X}^\kappa(\cC^{s}_{nsep}(\bar X)) = \sum_{\gamma \in \cP^s_{nsep}(X)} F(\ell_\gamma(X)). \eeq
    Applying Mirzakhani's integration formula (Theorem \ref{Mirzakhani's integration formula}) to $F$, we obtain
    \beq  \EgnL \Big[ \mu_{\bar X}^\kappa(\cC^{s}_{nsep}(\bar X)) \Big] = \int_0^\infty F(x) \frac{V_{g-1,n+2}(\mathbf{L},x,x)}{V_{g,n}(\mathbf{L})} xdx. \eeq
    By Lemma \ref{Volume bound II} and the dominated convergence theorem, it follows that
    \beq \lim_{g \to \infty} \int_0^\infty F(x) \frac{V_{g-1,n+2}(\mathbf{L},x,x)}{V_{g,n}(\mathbf{L})} xdx = \int_0^\infty F(x) \cdot \frac{(e^x - 1)^2}{x e^x} dx .\eeq
    By the definition of $F$, the RHS integral is precisely $c(k)$.
    \ep

    \begin{rem*} We now derive an expansion of \(c(k)\) in the same form as in Theorem \ref{Main theorem}. Starting from its integral representation,
    \beq c(k) = \sum_{m=1}^\infty \int_0^\infty \frac{e^{kx}}{m(e^{mx}-1)} \frac{(e^x-1)^2}{xe^x}dx,\ k < 0. \eeq
    Let $J_m(k)$ denote $m$ times the $m$-th term in the summation, namely,
    \beq J_m(k) = \int_0^\infty \frac{e^{kx}}{e^{mx} - 1} \frac{(e^x-1)^2}{xe^x}\ dx = \sum_{n=1}^\infty \int_0^\infty e^{-mnx} \cdot e^{kx} \frac{(e^x-1)^2}{xe^x} dx. \eeq
    Using the identity $\int_0^\infty \frac{e^{-ax} - e^{-bx}}{x}dx = \log \frac{b}{a}$, we have
    \beq \begin{split} J_m(k) &= \sum_{n=1}^\infty \int_0^\infty \frac{e^{-(mn-k-1)x} -2 e^{-(mn-k)x} + e^{-(mn-k+1)x}}{x}dx \\ &=  \sum_{n=1}^\infty \log \Big( \frac{(mn-k)^2}{(mn-k-1)(mn-k+1)} \Big). \end{split} \eeq
    Consequently,
    \beq \begin{split} c(k) & = \sum_{m,n=1}^\infty \frac{1}{m} \Big[2\log \Big(1 - \frac{k}{mn} \Big) - \log \Big(1 - \frac{k+1}{mn} \Big) - \log \Big(1 - \frac{k-1}{mn} \Big) \Big] \\ & = \sum_{N=1}^\infty \frac{\sigma(N)}{N} \Big[2\log \Big(1 - \frac{k}{N} \Big) - \log \Big(1 - \frac{k+1}{N} \Big) - \log \Big(1 - \frac{k-1}{N} \Big) \Big], \end{split} \eeq
    where $\sigma(N) = \sum_{d | N} d$ denotes the sum of divisors function. It's not hard to see that $\sigma(N) < N \log N$ for any $N \geq 2$. Assuming $k > -1$, we split the summation into the $N=1$ term and the $N \geq 2$ terms, then for each $N \geq 2$ term, expand the logarithm via its Taylor series to obtain
    \beq c(k) = \log \frac{(k-1)^2}{k(k-2)} -     \sum_{N=2}^\infty \sum_{j=2}^\infty \frac{2k^j - (k+1)^j -(k-1)^j}{j} \cdot \frac{\sigma(N)}{N^{j+1}}, \eeq
    where the double series converges absolutely for $-1 < k < 0$. Recall that
    \beq \sum_{N=1}^\infty \frac{\sigma(N)}{N^{j+1}} = \zeta(j+1) \zeta(j),\ \forall j \geq 2. \eeq
    Hence, by commuting the order of summation, we have
    \beq c(k) = \log \frac{(k-1)^2}{k(k-2)} - \sum_{j=2}^\infty \frac{2k^j - (k+1)^j -(k-1)^j}{j} \cdot \Big[ \zeta(j+1) \zeta(j) - 1 \Big], \eeq
    which is precisely the same expression as in Theorem \ref{Main theorem}.
    \end{rem*}

    \vone
    
    If $k=0$, following the same steps as in Proposition \ref{nsep, m=all, k<0} and the subsequent remark yield the following result:
    
    \bpro \label{nsep, m geq 2, k=0} Fix $n \geq 0$, and assume $|\mathbf{L}|^2 = o(g) $. Then
    \beq \lim_{g \to \infty} \frac{1}{V_{g,n}(\mathbf{L})} \int_{\cM_{g,n}(\mathbf{L})} \sum_{m=2}^\infty \sum_{\gamma \in \cP^s_{nsep}(X)} \frac{1}{m (e^{ m\ell_\gamma(X)} - 1)} dX = c',  \eeq
    where $c'$ is a constant given by
    \beqar c'&=&\int_0^\infty \sum_{m=2}^\infty \frac{1}{m(e^{mx} - 1)} \frac{(e^x - 1)^2}{x e^x} dx \\
    &=& \sum_{m=2}^\infty \frac{\zeta(2m)(\zeta(2m+1)-1)}{m}\approx 0.0234\cdots. \eeqar
    \epro

Now we estimate the leading term.

    \bpro \label{nsep, m = 1, k = 0} Fix $n \geq 0$ and assume that $|\mathbf{L}|^2 = o(g)$. Then
    \beq  \frac{1}{V_{g,n}(\mathbf{L})} \int_{\cM_{g,n}(\mathbf{L})} \sum_{\gamma \in \cP^s_{nsep}(X)} \frac{1}{e^{ \ell_\gamma(X)} - 1} dX \sim \frac{1}{2} \log g. \eeq
    Moreover, if $|\mathbf{L}|^2 = o(\frac{g}{\log g})$, then there is a constant $C_n $ depending only on $n$ such that for any sufficiently large $g$, 
    \beq  \Bigg| \frac{1}{V_{g,n}(\mathbf{L})} \int_{\cM_{g,n}(\mathbf{L})} \sum_{\gamma \in \cP^s_{nsep}(X)} \frac{1}{e^{ \ell_\gamma(X)} - 1} dX  - \frac{1}{2} \log g \Bigg| \leq C_n. \eeq
    \epro

    \bp Define the integration function
    \beq F(x) = \frac{1}{e^x - 1}. \eeq
    Applying Mirzakhani's integration formula (Theorem \ref{Mirzakhani's integration formula}) to $F$, we obtain
    \be \EgnL \Big[ \sum_{\gamma \in \cP^s_{nsep}(X)} F(\ell_\gamma(X)) \Big] = \int_0^\infty F(x) \frac{V_{g-1,n+2}(\mathbf{L},x,x)}{V_{g,n}(\mathbf{L})} xdx. \label{Integration formula nsep m = 1} \ene
    We now split the RHS of \eqref{Integration formula nsep m = 1} into two parts:
    \beq \mathrm{I} = \int_{0}^{\sqrt{\frac{g}{4c_n}}} F(x) \frac{V_{g-1,n+2}(\mathbf{L},x,x)}{V_{g,n}(\mathbf{L})} xdx, \eeq
    \beq \mathrm{II} = \int_{\sqrt{\frac{g}{4c_n}}}^\infty F(x) \frac{V_{g-1,n+2}(\mathbf{L},x,x)}{V_{g,n}(\mathbf{L})} xdx. \eeq
    For the second part, we apply part (2) of Lemma \ref{Volume bound II} to get the bound
    \beq \begin{split} \mathrm{II} &\leq \alpha^2 (2g+n-2)^2 \Big(1 - c_n \frac{|\mathbf{L}|^2 }{g} \Big)^{-1} \Big(1 + \frac{c_n}{g} \Big) \\ & \times \int_{\sqrt{\frac{g}{4c_n}}}^\infty \frac{F(x)}{x^4} \Big( \frac{\sinh(x/2)}{x/2} \Big)^2 \cdot xdx. \end{split} \eeq
    The integral can be evaluated explicitly:
    \beq \int_{\sqrt{\frac{g}{4c_n}}}^\infty \frac{F(x)}{x^4} \Big( \frac{\sinh(x/2)}{x/2} \Big)^2 \cdot xdx = \int_{\sqrt{\frac{g}{4c_n}}}^\infty \frac{1-e^{-x}}{x^5} dx \leq \frac{4 c_n^2}{g^2}.  \eeq
    Consequently, 
    \beq \mathrm{II} \leq 4 c_n^2 \alpha^2 \Big( \frac{2g+n-2}{g} \Big)^2 \Big(1 - c_n \frac{|\mathbf{L}|^2 }{g} \Big)^{-1} \Big(1 + \frac{c_n}{g} \Big) \prec 1.  \eeq
    For the first part, we compare $\mathrm{I}$ with
    \beq \mathrm{I'} = \int_0^{\sqrt{\frac{g}{4c_n}}} F(x) \Big( \frac{\sinh(x/2)}{x/2} \Big)^2 \cdot xdx = \int_0^{\sqrt{\frac{g}{4c_n}}} \frac{1-e^{-x}}{x} dx,  \eeq
    whose leading term is $\int_1^{\sqrt{\frac{g}{4c_n}}}\frac{1}{x}dx=\frac{1}{2} \log g-\frac{1}{2} \log (4c_n)$. It is straightforward to see that there exists a constant $c_n'>0$ depending only on $n$ such that for any $g \geq 2$,
    \beq \Big|\mathrm{I'} - \frac{1}{2} \log g \Big| \leq c_n' . \eeq
    By part (1) of Lemma \ref{Volume bound II}, we have
    \be \label{comparing I and I'} \begin{split}
    |\mathrm{I} - \mathrm{I'}| &\leq \frac{3c_n}{g} \int_0^{\sqrt{\frac{g}{4c_n}}} F(x) \Big( \frac{\sinh(x/2)}{x/2} \Big)^2 \cdot (1 + |\mathbf{L}|^2 + x^2) xdx \\ & = \frac{3c_n(1 + |\mathbf{L}|^2)}{g} \mathrm{I'} + \frac{3c_n}{g} \int_0^{\sqrt{\frac{g}{4c_n}}} x(1-e^{-x}) dx \\ & \leq \frac{3c_n(1 + |\mathbf{L}|^2)}{g}\mathrm{I'} + \frac{3}{8}. \end{split} \ene
Since $c_n(1 + |\mathbf{L}|^2) = o(g)$, it follows from \eqref{comparing I and I'} that
 \beq \lim_{g \to \infty} \frac{\mathrm{I}}{\mathrm{I'}} = 1 , \eeq
    and therefore $\mathrm{I} \sim \frac{1}{2} \log g$, so $\mathrm{I} + \mathrm{II} \sim \frac{1}{2} \log g$, i.e.
    \beq \frac{1}{V_{g,n}(\mathbf{L})} \int_{\cM_{g,n}(\mathbf{L})} \sum_{\gamma \in \cP^s_{nsep}(X)} \frac{1}{e^{ \ell_\gamma(X)} - 1} dX \sim \frac{1}{2} \log g. \eeq
    If, in addition, $|\mathbf{L}|^2 = o(\frac{g}{\log g})$, then $\frac{3c_n(1 + |\mathbf{L}|^2)}{g} \mathrm{I'} \to 0$ as $g\to \infty$. Thus, there is a constant $C_n$ that depends only on $n$, such that for all sufficiently large $g$, we have
    \beq \Bigg| \frac{1}{V_{g,n}(\mathbf{L})} \int_{\cM_{g,n}(\mathbf{L})} \sum_{\gamma \in \cP^s_{sep}(X)} \frac{1}{e^{ \ell_\gamma(X)} - 1} dX  - \frac{1}{2} \log g \Bigg| \leq C_n. \eeq
    This completes the proof.
    \ep    

    \vone

    \subsection{Proofs of Theorems \ref{main theorem simple curves} and \ref{Total mass of simple curves, case kappa = 0, duplicate}} 
    
    \bp[Proof of Theorem \ref{main theorem simple curves}] By Propositions \ref{sep, m = all, k = 0} and \ref{nsep, m=all, k<0}, for $\kappa >  0$ and $|\mathbf{L}|^2 = o(g)$, we have
    \beq \lim_{n \to \infty} \EgnL \Big[\mu_{\bar X}^\kappa(\cC^s_{sep}(\bar X)) \Big]  = 0,  \eeq
    and
    \beq \lim_{n \to \infty} \EgnL \Big[\mu_{\bar X}^\kappa(\cC^s_{nsep}(\bar X)) \Big] =  c \Bigg(\frac{1}{2} - \sqrt{\frac{1}{4} + \kappa} \Bigg). \eeq
    Since $\mu_{\bar X}^\kappa(\cC^s(\bar X)) = \mu_{\bar X}^\kappa(\cC^s_{sep}(\bar X)) + \mu_{\bar X}^\kappa(\cC^s_{nsep}(\bar X))$, combining the two limits above completes the proof. 
    \ep

   \bp[Proof of Theorem \ref{Total mass of simple curves, case kappa = 0, duplicate}]
    For $\kappa = 0$ and $|\mathbf{L}|^2 = o(g)$, Proposition \ref{sep, m = all, k = 0} still gives 
    \beq \lim_{n \to \infty} \EgnL \Big[\mu_{\bar X}(\cC^s_{sep}(\bar X)) \Big] = 0,  \eeq
    and by Propositions \ref{nsep, m geq 2, k=0} and \ref{nsep, m = 1, k = 0}, 
    \be \lim_{g \to \infty} \frac{\EgnL \Big[\mu_{\bar X}(\cC^s_{nsep}(\bar X)) \Big]}{\frac{1}{2}\log g} = 1. \ene
    Moreover, if $|\mathbf{L}|^2 = o(\frac{g}{\log g})$, then Proposition \ref{nsep, m = 1, k = 0} implies 
    \be \Bigg| \EgnL \Big[\mu_{\bar X}(\cC^s(\bar X)) \Big]  - \frac{1}{2} \log g \Bigg| \prec 1.  \ene
    Putting all these together completes the proof.
    \ep

    \vone

    \section{The total mass of $\mu_{\bar X}^\kappa$ over non-simple closed curves}\label{sec4-nscc}

    Our main result for this section is as follows:

    \bt \label{main theorem nonsimple curves}
    Fix $\kappa > 0$ and $n \geq 0$. Let $\mathbf{L} =\mathbf{L}(g) \in \R_{\geq 0}^n$  satisfy $|\mathbf{L}|^2 = o(g)$. Then for any $0 < \eps < 1$, as $g \to \infty$,
    \be \EgnL \Big[\mu_{\bar X}^\kappa(\cC^{ns}(\bar X)) \Big] \prec \frac{1}{g^{1 - \eps}}. \ene
    \et

    Now we prove Theorem \ref{Main theorem}, assuming Theorem \ref{main theorem nonsimple curves}.
    \bp[Proof of Theorem \ref{Main theorem}]
    The conclusion follows clearly from Theorems \ref{main theorem simple curves} and \ref{main theorem nonsimple curves}.
    \ep

    Before showing Theorem \ref{main theorem nonsimple curves}, we make the following preparations.
    
    Let $Y_0 \subseteq X$ be a geodesic subsurface that satisfies the following conditions: \ben
    \item $Y_0 \cong S_{g_0,k}$ for some $g_0 ,k \geq 0$ and $m = |\chi(Y_0)| = 2g_0 + k -2 \geq 1$.
    \item The boundary $\p Y_0$ contains $n_0 \geq 0$ pairs of closed geodesics, where $2n_0 \leq k$. Each such pair corresponds to a single simple non-separating closed geodesic in $X$.
    \item $Y_0$ has $n_f \geq 0$ hyperbolic ends that are geodesic boundaries of $X$ or cusps of $X$. Denote by $L_{i_1},\cdots,L_{i_{n_f}} \geq 0$ the geodesic lengths of these ends, where $L_{i_j} = 0$ means that the $i_j$-th end is a cusp of $X$. Let $\mathbf{L}_{n_f} = (L_{i_1},\cdots,L_{i_{n_f}}) \in \R_{\geq 0}^{n_f}$, it is clear that
    \beq  |\mathbf{L}_{n_f}|_1 =\sum_{j=1}^{n_f}L_{i_j}= \ell_X(\p Y_0 \cap \p X). \eeq
    \item The complement decomposes as $X \backslash Y_0 \cong \sqcup_{i=1}^q S_i$, where $S_i \cong S_{g_i,n_i}$. In particular, 
    \beq n = \sum_{i=1}^q n_i + 2n_f + 2n_0 - k,\ \sum_{i=1}^q (2g_i - 2 + n_i) = 2g + n - 2 - m. \eeq
    \item For each $1 \leq i \leq q$, $S_i$ has $n_{i,f}$ hyperbolic ends that are geodesic boundaries of $X$ or cusps of $X$, where $0 \leq n_{i,f} \leq n_i$. Let $\mathbf{L}_{n_i,f} \in \R_{\geq 0}^{n_{i,f}}$ be the geodesic lengths of these ends, then $|\mathbf{L}_{n_i,f}|_1 = \ell_X(\p S_i \cap \p X)$, and
    \beq \sum_{i=1}^q n_{i,f} + n_f = n. \eeq
    \een
    The following is a direct consequence of Mirzakhani's integration formula.

    \bpro Let $F \in C([0,2T])$ be a continuous function supported on $[0,2T]$, and let $Y_0 \subseteq X$ be as above. If $T < (g-1)\pi/2$, then
    \be \label{Over a fixed orbit} \begin{split}
        & \quad \quad \quad \int_{\cM_{g,n}(\mathbf{L})} \sum_{Y \in \mathrm{Mod}_{g,n} \cdot Y_0} F(\ell_X (\p Y))\ dX 
        \\ &= \frac{C(\p Y_0)}{|\mathrm{Sym}(Y_0)|}\int_{\R^{k - n_0 - n_f}_{\geq 0}} F(|\mathbf{L}_{n_f}|_1 + \sum_{i=1}^{n_0} 2x_i' + \sum_{i=1}^{q} \sum_{j=1}^{n_i - n_{i,f}} x_{i,j}) \\ & \times  V_{g_0,k}(x^{(0)},\mathbf{L}_{n_f}) V_{g_1,n_1}(x^{(1)},\mathbf{L}_{n_{1,f}}) \cdots V_{g_q,n_q}(x^{(q)},\mathbf{L}_{n_q,f}) 
        \\ & \times x_1' \cdots x_{n_0}' x_{1,1} \cdots x_{q,n_q - n_{q,f}} dx_1' \cdots dx_{n_0}' dx_{1,1} \cdots dx_{q,n_q - n_{q,f}}, 
    \end{split} \ene
    where $\mathrm{Sym}(Y_0)$ is the symmetry group of $\p Y_0$; $C(\p Y_0) \in (0,1]$ is a constant depending only on the topological type of $\p Y_0$, and 
    \beq V_{g_0,k}(x^{(0)},\mathbf{L}_{n_f}) V_{g_1,n_1}(x^{(1)},\mathbf{L}_{n_{1,f}}) \cdots V_{g_q,n_q}(x^{(q)},\mathbf{L}_{n_q,f}) =  V_{g,n}(\p Y_0 \backslash \p X,x) \eeq
    is the volume of the moduli space of Riemann surfaces homotopic to $X \backslash \p Y_0$, with boundary vector length given by $x \in \R^{k - n_0 - n_f}_{\geq 0}$.
    \epro

    \bp Recall that $\mathrm{Sub}_T(X)$ is the set of all subsurfaces $Y \subseteq X$ with geodesic boundary that satisfies $\area(Y) \leq 4T$ and $\ell_X(\p Y) \leq 2T$. If $2T < \pi(g-1)$, then the map $Y \mapsto \p Y$ is injective on $\mathrm{Sub}_T(X)$ (see \cite[Section 7.1]{Wu_2022}). It follows that
    \beqar \quad && \quad \int_{\cM_{g,n}(\mathbf{L})} \sum_{Y \in \mathrm{Mod}_{g,n} \cdot Y_0} F(\ell_X(\p Y))\ dX 
    \\ && = \int_{\cM_{g,n}(\mathbf{L})} \sum_{\Gamma \in \mathrm{Mod}_{g,n} \cdot \p Y_0} F(\ell_X(\Gamma))\ dX
    \\ && = \int_{\cM_{g,n}(\mathbf{L})} \sum_{\Gamma \in \mathrm{Mod}_{g,n} \cdot (\p Y_0 \backslash \p X)} F \big( |\mathbf{L}_{n_f}|_1 + |\ell_X(\Gamma)| \big)\ dX. \eeqar
    The conclusion then follows from Mirzakhani's integration formula (Theorem \ref{Mirzakhani's integration formula}).
    \ep

    For $T>0$, define the domain $\Delta \subseteq \R_{+}^{k - n_0 - n_f}$ via
    \be \label{The integration domain} \Delta = \Big\{(x_i',x_{i,j}) \in \R^{k - n_0 - n_f}_{\geq 0}: \sum_{i=1}^{n_0} 2x_i' + \sum_{i=1}^{q} \sum_{j=1}^{n_i - n_{i,f}} x_{i,j} \leq 2T \Big\}, \ene
    and equip $\R_+^{k - n_0 - n_f}$ with the measure
    \beq d\nu(x) = x_1' \cdots x_{n_0}' x_{1,1} \cdots x_{q,n_q - n_{q,f}} dx_1' \cdots dx_{n_0}' dx_{1,1} \cdots dx_{q,n_q - n_{q,f}}. \eeq

  \noindent It is easy to check
    \beq \nu(\Delta) = \frac{(2T)^{2(k - n_0 - n_f)}}{(2(k - n_0 - n_f))!}. \eeq 

    \bpro \label{Over a fixed topological type} 
    For any $\lambda > 0$, integers $g_0$ and $k$ satisfying $2g_0 + k -2 \geq 1$, and any $T < (g-1)\pi/2$, then we have
\beqar
&&  \quad \quad \quad \int_{\cM_{g,n}(\mathbf{L})} \sum_{
    \begin{subarray}{1}
        Y \in \mathrm{Sub}_T(X), \\ \ \  Y \simeq S_{g_0,k}
    \end{subarray}}  e^{- \lambda \ell_X(\p Y)}\ dX \\
 &&    \leq e^{(2-2\lambda)_+ T} \cdot \prod_{i=1}^n \frac{\sinh(L_i/2)}{L_i/2} \cdot \sum_{\ast} \frac{V_{g_0,k} \prod_{i=1}^q V_{g_i,n_i} }{n_0! \prod_{i=1}^q (n_i - n_{i,f})!}\cdot \nu(\Delta).
\eeqar
 Here, $\sum\limits_{\ast}$ runs over all $\mathrm{Mod}_{g,n}$-orbits of geodesic subsurfaces $Y_0 \cong S_{g_0,k} \subseteq X$, and we write $(2-2\lambda)_+=\max\left\{2-2\lambda, 0 \right\}$.
    \epro
    
    \bp The $n_0$ pairs of geodesics of $Y_0$ contribute a factor of $n_0!$ to $|\mathrm{Sym}(Y_0)|$, and the permutations of the $n_i - n_{i,f}$ boundary geodesics of $\p S_{g_i,n_i}$ lying in the interior of $X$ contribute a factor $(n_i - n_{i,f})!$ to $|\mathrm{Sym}(Y_0)|$, for each $1 \leq i \leq q$. This yields the lower bound 
    \beq |\mathrm{Sym}(Y_0)| \geq n_0! \cdot \prod_{i=1}^q (n_i - n_{i,f})!. \eeq
    By the Weil-Petersson volume bounds \eqref{Mirzakhani's volume comparison formula} and \eqref{NWX volume comparison}, we deduce
    \beq \frac{V_{g,n}(\p Y_0 \backslash \p X,x)}{V_{g_0,k} \prod_{i=1}^q V_{g_i,n_i}} \leq \prod_{i=1}^n \frac{\sinh(L_i/2)}{L_i/2} \cdot \exp \Big(\sum_{i=1}^{n_0} x_i' + \sum_{i=1}^{q} \sum_{j=1}^{n_i - n_{i,f}} x_{i,j}\Big). \eeq
    Set $F(x) = e^{- \lambda x} \cdot \mathbf{1}_{[0,2T]}(x)$. Substituting it into \eqref{Over a fixed orbit} yields
    \beq \begin{split}
        & \quad \quad \int_{\cM_{g,n}(\mathbf{L})} \sum_{ Y \in \mathrm{Mod}_{g,n} \cdot Y_0} e^{-\lambda \ell_X(\p Y)} \mathbf{1}_{ \{\ell_X(\p Y) \leq 2T\} }\  dX
        \\ &\leq \frac{1}{n_0! \prod_{i=1}^q(n_i - n_{i,f})!} \prod_{i=1}^n \frac{\sinh(L_i/2)}{L_i/2} \cdot V_{g_0,k} \prod_{i=1}^q V_{g_i,n_i} \\ & \times \int_{\Delta\subset \R^{k - n_0 - n_f}_{\geq 0}} \exp \Big((1 -2 \lambda) \sum_{i=1}^{n_0} x_i' + (1 - \lambda) \sum_{i=1}^{q} \sum_{j=1}^{n_i - n_{i,f}} x_{i,j} \Big) d\nu(x). 
    \end{split} \eeq
    Since $\supp(F) \subseteq [0,2T]$, the integral restricts to the domain $\Delta$. On this domain,
    \beq (1-2\lambda)\sum_{i=1}^{n_0} x_i' + (1-\lambda)\sum_{i=1}^{q} \sum_{j=1}^{n_i - n_{i,f}} x_{i,j} \leq (2 - 2\lambda)_+ T. \eeq
    Bounding the integral over $\Delta$ by $e^{(2-2\lambda)_+ T} \nu(\Delta)$, we obtain
    \beqar 
    && \quad \quad \int_{\cM_{g,n}(\mathbf{L})} \sum_{
    \begin{subarray}{1} Y \in \mathrm{Mod}_{g,n} \cdot Y_0, \\ \   Y \in \mathrm{Sub}_T(X)            
    \end{subarray}} e^{- \lambda \ell_X(\p Y)}\ dX\\
    &&\leq e^{(2-2\lambda)_+ T} \cdot \prod_{i=1}^n \frac{\sinh(L_i/2)}{L_i/2} \cdot \frac{V_{g_0,k} \prod_{i=1}^q V_{g_i,n_i} }{n_0! \prod_{i=1}^q (n_i - n_{i,f})!}\cdot \nu(\Delta). 
    \eeqar
    Summing over all $\mathrm{Mod}_{g,n}$-orbits of such subsurfaces completes the proof.
    \ep
    
    For the remainder of the argument, fix a constant $c_T > 0$ and set $T =
     c_T \cdot \log g$.

    \bpro \label{large Euler characteristic over a fixed type} Fix $\lambda > 0$ and nonnegative integers $n,g_0,k\geq 0$. There exist two universal constants $c,D > 0$ such that for any $g \geq \max(D,c_T^2)$ and every $\mathbf{L}\in \R^n_{\geq 0}$,
\beqar
&& \quad  \quad  \int_{\cM_{g,n}(\mathbf{L})} \sum_{
    \begin{subarray}{1}
        \ Y \in \mathrm{Sub}_T(X), \\ \ \ \  Y \simeq S_{g_0,k}
    \end{subarray}} e^{- \lambda \ell_X(\p Y)}\ dX \\
&&\leq  c \cdot n! e^{\max \left\{4-2\lambda,2 \right\} T} \prod_{i=1}^n \frac{\sinh(L_i/2)}{L_i/2} \cdot  V_{g_0, k} W_{2g + n - 2g_0 - k}.
\eeqar
    \epro

    \bp The sum $\sum\limits_{\ast}$ in Proposition 4.3 runs over all $\mathrm{Mod}_{g,n}$-orbits of embedded subsurfaces $S_{g_0,k} \hookrightarrow X$.
    Each such orbit fixes a collection of parameters \(q,n_f,n_0,g_i,n_i,n_{i,f}\) with \(1\le i\le q\), subject to the two topological identities
    \be \label{constraint} n = \sum_{i=1}^q n_i + 2n_f + 2n_0 - k = \sum_{i=1}^q n_{i,f} + n_f, \ene
    \be \sum_{i=1}^q (2g_i - 2 + n_i) = 2g + n - 2g_0 - k. \ene
    For given $g_0,k,q, n_f,n_0,g_i,n_i,n_{i,f}$, the total number of distinct $\mathrm{Mod}_{g,n}$ orbits is bounded from above by
    \beq  N  = \bpmatrix n \\ n_f \epmatrix \bpmatrix n - n_f \\ n_{1,f} \epmatrix \cdots \bpmatrix n_{q-1,f} + n_{q,f} \\ n_{q-1,f} \epmatrix \bpmatrix n_{q,f} \\ n_{q,f} \epmatrix = \frac{n!}{n_f! \prod_{i=1}^q n_{i,f}!}.  \eeq
    Using this count, we expand $(\ast)$ as a sum over the above parameters:
    \beq \label{summation over parameters} \sum_{\ast}  \frac{V_{g_0,k} \prod_{i=1}^q V_{g_i,n_i}}{n_0! \prod_{i=1}^q (n_i - n_{i,f})!} = \sum_q \sum_{n_i,g_i} \sum_{n_0,n_f,n_{i,f}} \frac{n! \cdot V_{g_0,k} \prod_{i=1}^q V_{g_i,n_i}}{n_0! n_f! \prod_{i = 1}^q n_{i,f}! (n_i - n_{i,f})!}. \eeq
    For fixed $q$ and $n_1,\cdots,n_q$, estimate \eqref{NWX volume summation estimate} gives
    \beq \sum_{\{g_i\}} V_{g_1,n_1} \cdots V_{g_q,n_q} \leq c \Big(\frac{D}{2g+n-2g_0 - k}\Big)^{q-1} W_{2g+n-2g_0-k}. \eeq
    Since $\area(Y) = 2\pi|\chi(Y)| \leq 4T$, we deduce a uniform bound on $k$ and $g_0$:
    \beq 2g_0 + k \leq |\chi(Y)| + 2 \leq \frac{2T}{\pi} + 2. \eeq
    After possibly enlarging the constant $D$, we have \(2g_0+k\le g\) for all \(g\ge\max(D,c_T^2)\), which further yields
    \beq \sum_{\{g_i\}} V_{g_1,n_1} \cdots V_{g_q,n_q} \leq c \Big(\frac{D}{g}\Big)^{q-1} W_{2g + n - 2g_0 - k}. \eeq
    Consequently, the sum $\sum\limits_\ast$ is bounded from above by
    \beq  \sum_q \sum_{n_0,n_f} \sum_{n_i} \sum_{n_{i,f}} \frac{c \cdot n! V_{g_0,k} W_{2g + n - 2g_0 - k}}{n_0! n_f! \prod_{i = 1}^q n_{i,f}! (n_i - n_{i,f})!} \Big(\frac{D}{g}\Big)^{q-1},\eeq
    where all indices satisfy the sole constraint \eqref{constraint}. Direct combinatorial estimation gives
    \beq \sum_{n_0,n_f} \sum_{n_1,\cdots,n_q} \sum_{n_{i,f}} \frac{1}{n_0! n_f! \prod_{i = 1}^q n_{i,f}! (n_i - n_{i,f})!} \leq e^{2 + 2q}. \eeq
  As a consequence,
    \beq \sum_{\ast}  \frac{V_{g_0,k} \prod_{i=1}^q V_{g_i,n_i}}{n_0! \prod_{i \geq 1} (n_i - n_{i,f})!} \leq c \cdot e^4  n! \cdot V_{g_0,k} W_{2g + n - 2g_0 - k} \sum_{q} \Big(\frac{e^2 D}{g}\Big)^{q-1}. \eeq
    The last series converges whenever $g \geq e^2 \cdot D$. Hence, for sufficiently large $g$, after enlarging the constant $c$ by a constant multiple, we have
    \beq \sum_{\ast}  \frac{V_{g_0,k} \prod_{i=1}^q V_{g_i,n_i}}{n_0! \prod_{i \geq 1} (n_i - n_{i,f})!} \leq c \cdot n! \cdot V_{g_0,k} W_{2g + n - 2g_0 - k}. \eeq 
    Finally, combining the bound
    \beq \nu(\Delta) = \frac{(2T)^{2(k - n_0 - n_f)}}{(2(k - n_0 - n_f))!} \leq e^{2T} \eeq
    with the estimates derived in Proposition \ref{Over a fixed topological type} yields
    \beqar
    && \quad \quad \int_{\cM_{g,n}(\mathbf{L})} \sum_{
    \begin{subarray}{1}
        \ Y \in \mathrm{Sub}_T(X), \\ \ \ \  Y \simeq S_{g_0,k}
    \end{subarray}} e^{- \lambda \ell_X(\p Y)}\ dX \\
    && \leq c \cdot n! e^{\max \left\{4-2\lambda,2 \right\} T} \prod_{i=1}^n \frac{\sinh(L_i/2)}{L_i/2} \cdot  V_{g_0, k} W_{2g + n - 2g_0 - k}. 
    \eeqar
    This completes the proof.
    \ep

    \subsection{Two upper bounds}
    In this subsection, we prove the following two upper bounds, which will be invoked in subsequent arguments.
    \bpro \label{Bounding exponential length function} Fix $A \geq 1, n \geq 0$ and $\lambda > 0$. Let $c_T > 0$ and define $T = c_T \log g$. If $|\mathbf{L}|^2 = o(g)$ as $g \to \infty$, then the following two asymptotic bounds hold:
    \be \label{first bound for section 4.1} \EgnL\Big[\sum_{
    \begin{subarray}{1}
        \ Y \in \mathrm{Sub}_T(X), \\ \ |\chi(Y)| \geq A + 1
    \end{subarray}} e^{-\lambda \ell_X(\p Y)} \Big]  \prec \frac{T e^{\max \left\{4-2\lambda,2 \right\} T}}{g^{A+1}}, \ene
    and
    \be \EgnL\Big[\sum_{
    \begin{subarray}{1}
        \ Y \in \mathrm{Sub}_T(X), \\ \ 1 \leq |\chi(Y)| \leq A
    \end{subarray}} e^{-\lambda \ell_X(\p Y)} \Big]  \prec \frac{T^{4A+2} e^{(1 - 2 \lambda)_+ T}}{g}. \ene
    \epro
    We first derive an estimate for subsurfaces with relatively small (absolute values of) Euler characteristics:

    \bpro \label{Small Euler characteristics over a fixed type} Fix $\lambda > 0$ and $n \geq 0$. For non-negative integers $g_0, k \geq 0$, there exists a constant $c(m.n)>0$ depending on $m = 2g_0 + k -2$ and $n$, such that for all sufficiently large $g$, 
\beqar
&& \quad \quad \int_{\cM_{g,n}(\mathbf{L})} \sum_{
    \begin{subarray}{1}
        \ Y \in \mathrm{Sub}_T(X), \\ \ \ \  Y \simeq S_{g_0,k}
    \end{subarray}} e^{- \lambda \ell_X(\p Y)}\ dX\\
    &&\leq c(m,n) \cdot T^{4m+2} e^{(1-2\lambda)_+ T} \cdot \frac{V_{g,n}}{g^m} \prod_{i=1}^n \frac{\sinh(L_i/2)}{L_i/2}. 
\eeqar
    \epro

    \bp From the proof of \eqref{Over a fixed orbit}, we have the estimate
    \be \label{Over a fixed orbit 2}\begin{split}
        & \quad \quad \int_{\cM_{g,n}(\mathbf{L})} \sum_{Y \in \mathrm{Mod}_{g,n} \cdot Y_0} e^{-\lambda \ell_X (\p Y)} \mathbf{1}_{ \{ \ell_X(\p Y) \leq 2T \}}\ dX 
        \\ &\leq \frac{1}{n_0! \prod_{i=1}^q (n_i - n_{i,f})!} \int_{\Delta} \exp \Big(-\lambda \sum_{i=1}^{n_0} 2x_i' -\lambda \sum_{i=1}^{q} \sum_{j=1}^{n_i - n_{i,f}} x_{i,j} \Big) \\ & \times  V_{g_0,k}(x^{(0)},\mathbf{L}_{n_f}) V_{g_1,n_1}(x^{(1)},\mathbf{L}_{n_{1,f}}) \cdots V_{g_q,n_q}(x^{(q)},\mathbf{L}_{n_q,f}) 
        \\ & \times x_1' \cdots x_{n_0}' x_{1,1} \cdots x_{q,n_q - n_{q,f}} dx_1' \cdots dx_{n_0}' dx_{1,1} \cdots dx_{q,n_q - n_{q,f}}. 
    \end{split} \ene
    By \cite[Theorem 1.1]{Mirzakhani:2007a} we know that 
    $V_{g_0,k}(2y_1,\cdots,2y_k)$ is a polynomial in $y_1^2,\cdots,y_k^2$ of degree $3g_0 + k - 3 = \frac{3m-k}{2}$. Whose coefficient of $y_1^{2d_1} \cdots y_k^{2d_k}$ is bounded from above by $\frac{V_{g_0,k}}{(2d_1+1)! \cdots (2d_k+1)!}$ (see e.g. \cite[Remark on Page 286]{Mir13a}). Then it is not hard to deduce that (see e.g. \cite[Lemma 22]{Nie_2023}) there exists a constant $c(m)>0$ depending only on $m$ such that 
    \beq \begin{split} V_{g_0,k}(x^{(0)},\mathbf{L}_{n_f}) &\leq c(m) \cdot (1 + ||x^{(0)}||_\infty)^{3m-k} \cdot \sum_{\mathbf{d} \in \N^{n_f}} \prod_{j=1}^{n_f} \frac{L_{i_j}^{2d_j}}{(2d_j + 1)!}
    \\ & = c(m)\prod_{j=1}^{n_f} \frac{\sinh(L_{i_j}/2)}{L_{i_j}/2} \cdot (1 + ||x^{(0)}||_\infty)^{3m-k}. \end{split} \eeq
    For any $x \in \Delta$, the definition \eqref{The integration domain} of $\Delta$ enforces  $||x||_\infty \leq 2T$. Since $g > e^{1/c_T}$ implies $T >1$; after possibly enlarging the constant \(c(m)\), we simplify to 
    \beq V_{g_0,k}(x^{(0)},\mathbf{L}_{n_f}) \leq c(m) \prod_{j=1}^{n_f} \frac{\sinh(L_{i_j}/2)}{L_{i_j}/2} \cdot T^{3m-k}. \eeq
    Next we apply estimate \eqref{Mirzakhani's volume comparison formula} to each volume term $V_{g_i,n_i}(x^{(i)},\mathbf{L}_{n_{i,f}})$ for $1 \leq i \leq q$, yielding
    \beq \frac{V_{g,n}(\p Y_0 \backslash \p X,x)}{\prod_{i=1}^q V_{g_i,n_i}} \leq c(m) T^{3m-k} \prod_{i=1}^n \frac{\sinh(L_i/2)}{L_i/2} \cdot \exp \Big(\frac{1}{2} \sum_{i=1}^{q} \sum_{j=1}^{n_i - n_{i,f}} x_{i,j}\Big). \eeq
    Substituting this bound back into \eqref{Over a fixed orbit 2} gives
    \beq \begin{split}
        & \quad \quad \int_{\cM_{g,n}(\mathbf{L})}\sum_{ Y \in \mathrm{Mod}_{g,n} \cdot Y_0} e^{-\lambda \ell_X(\p Y)} \mathbf{1}_{ \{\ell_X(\p Y) \leq 2T\} }\  dX
        \\ &\leq \frac{c(m) T^{3m-k}}{n_0! \prod_{i=1}^q(n_i - n_{i,f})!} \prod_{i=1}^n \frac{\sinh(L_i/2)}{L_i/2} \cdot \prod_{i=1}^q V_{g_i,n_i}
        \\ & \times \int_{\Delta} \exp \Big(-2 \lambda \sum_{i=1}^{n_0} x_i' - \frac{1-2\lambda}{2} \sum_{i=1}^{q} \sum_{j=1}^{n_i - n_{i,f}} x_{i,j} \Big) d\nu(x). 
    \end{split} \eeq
    Since $x \in \Delta$, we have
    \beq -2\lambda \sum_{i=1}^{n_0} x_i' + \frac{1 - 2\lambda}{2} \sum_{i=1}^{q} \sum_{j=1}^{n_i - n_{i,f}} x_{i,j} \leq (1-2\lambda)_+ T. \eeq
    Bounding the integral over $\Delta$ by $e^{(1-2\lambda)_+ T} \cdot \nu(\Delta)$ leads to
    \beqar 
    && \quad \quad \int_{\cM_{g,n}(\mathbf{L})} \sum_{
    \begin{subarray}{1} Y \in \mathrm{Mod}_{g,n} \cdot Y_0, \\ \   Y \in \mathrm{Sub}_T(X)            
    \end{subarray}} e^{- \lambda \ell_X(\p Y)}\ dX \\
    &&\leq c(m) T^{3m-k} e^{(1-2\lambda)_+ T} \prod_{i=1}^n \frac{\sinh(L_i/2)}{L_i/2} \cdot \frac{\prod_{i=1}^q V_{g_i,n_i}  }{n_0! \prod_{i=1}^q (n_i - n_{i,f})!}\cdot \nu(\Delta). 
    \eeqar
    We next sum over all orbits. As before, this summation expands as a sum over the parameters $q,n_f,n_0,g_i,n_i,n_{i,f}$:
    \beq \sum_{\ast}  \frac{\prod_{i=1}^q V_{g_i,n_i}}{n_0! \prod_{i=1}^q (n_i - n_{i,f})!} = \sum_q \sum_{n_i,g_i} \sum_{n_0,n_f,n_{i,f}} \frac{n! \cdot \prod_{i=1}^q V_{g_i,n_i}}{n_0! n_f! \prod_{i = 1}^q n_{i,f}! (n_i - n_{i,f})!}. \eeq
    By the Weil-Petersson volume bound \eqref{Euler char lowers by m}, for any fixed $q$ and $\{n_i\}_{i=1}^q$, whenever $g$ is sufficiently large, we have
    \beq \sum_{\{g_i\}} \prod_{i=1}^q V_{g_i,n_i} \leq c(m,n) \frac{V_{g,n}}{g^m}. \eeq
    Repeating the combinatorial counting argument used in the proof of Proposition \ref{large Euler characteristic over a fixed type} leads to
    \beq \sum_{\ast}  \frac{\prod_{i=1}^q V_{g_i,n_i}}{n_0! \prod_{i \geq 1} (n_i - n_{i,f})!} \leq c(m,n) \cdot n! \cdot \frac{V_{g,n}}{g^m} \sum_{q} (e^2 D)^{q-1}. \eeq
    With a slight abuse of notation, we absorb the factor $n!$ into the constant $c(m,n)$. Recall \(m=2g_0+k-2\), which implies \(k\le m+2\) and in turn \(q\le k\le m+2\). Thus,
    \beq \sum_{\ast}  \frac{\prod_{i=1}^q V_{g_i,n_i}}{n_0! \prod_{i \geq 1} (n_i - n_{i,f})!} \leq c(m,n) \cdot \frac{V_{g,n}}{g^m}. \eeq
    We now substitute the volume bound for the domain measure
    \beq \nu(\Delta) = \frac{(2T)^{2(k - n_0 - n_f)}}{(2(k - n_0 - n_f))!} \leq T^{2k}. \eeq
    Combining all preceding estimates, we conclude
\beqar
&& \quad \quad     \int_{\cM_{g,n}(\mathbf{L})} \sum_{
    \begin{subarray}{1}
        \ Y \in \mathrm{Sub}_T(X), \\ \ \  Y \simeq S_{g_0,k}
    \end{subarray}} e^{- \lambda \ell_X(\p Y)}\ dX \\
    &&\leq c(m,n) \cdot T^{4m+2} e^{(1-2\lambda)_+ T} \cdot \frac{V_{g,n}}{g^m} \prod_{i=1}^n \frac{\sinh(L_i/2)}{L_i/2}. 
    \eeqar
    This completes the proof.
    \ep

    Now we are ready to prove Proposition \ref{Bounding exponential length function}.
    
    \bp[Proof of Proposition \ref{Bounding exponential length function}] By \cite[Lemma 23 (2)]{Nie_2023}, there is a constant $c(A)>0$ depending only on $A$ such that for all sufficiently large $g$,
    \be \label{NWX Lemma 23 (2)} \sum_{m = A}^{[2T/\pi]} W_m W_{2g+n-2-m} \leq c(A) \cdot \frac{W_{2g+n-2}}{g^A}. \ene
    By \cite[Lemma 23 (1)]{Nie_2023}, there is a universal constant $c>0$ such that for any relevant $g_0,k$, $V_{g_0,k} \leq c W_{2g_0 + k -2}$. This together with Proposition \ref{large Euler characteristic over a fixed type} and inequality \eqref{NWX Lemma 23 (2)} implies
    \beqar && \quad \quad \quad \int_{\cM_{g,n}(\mathbf{L})} \sum_{
    \begin{subarray}{1}
        \ Y \in \mathrm{Sub}_T(X), \\ \ |\chi(Y)| \geq A + 1
    \end{subarray}} e^{- \lambda \ell_X(\p Y)}\ dX 
    \\ && \prec \sum_{m=A+1}^{[2T/\pi]} m \cdot e^{\max \left\{4-2\lambda,2 \right\} T} \prod_{i=1}^n \frac{\sinh(L_i/2)}{L_i/2}  \cdot W_m W_{2g + n + 2 - m}
    \\ && \prec T e^{\max \left\{4-2\lambda,2 \right\} T} \frac{W_{2g+n-2}}{g^{A+1}} \prod_{i=1}^n \frac{\sinh(L_i/2)}{L_i/2}. \eeqar
    Since $|\mathbf{L}|^2 = o(g)$, the volume bounds \eqref{Volume asymptotics for adjacent g,n} and \eqref{NWX volume comparison} further imply
    \beq \int_{\cM_{g,n}(\mathbf{L})} \sum_{
    \begin{subarray}{1}
        \ Y \in \mathrm{Sub}_T(X), \\ \ |\chi(Y)| \geq A + 1
    \end{subarray}}  e^{- \lambda \ell_X(\p Y)}\ dX \prec T e^{\max \left\{4-2\lambda,2 \right\} T} \frac{V_{g,n}(\mathbf{L})}{g^{A+1}}. \eeq
    This completes the proof of the first asymptotic bound \eqref{first bound for section 4.1}.
    
    We next treat the regime \(1\le|\chi(Y)|\le A\). For every $m \geq 1$, Proposition \ref{Small Euler characteristics over a fixed type} together with the volume bound \eqref{NWX volume comparison} imply the upper estimate 
    \be \int_{\cM_{g,n}(\mathbf{L})} \sum_{
    \begin{subarray}{1}
        \ Y \in \mathrm{Sub}_T(X), \\ \ \  |\chi(Y)| = m
    \end{subarray}} e^{- \lambda \ell_X(\p Y)}\ dX \prec  T^{4m+2} e^{(1-2\lambda)_+ T} \frac{V_{g,n}(\mathbf{L})}{g^m}. \label{small Euler characteristics} \ene
    Summing \eqref{small Euler characteristics} over $1 \leq m \leq A$, we deduce
    \beq \begin{split} \int_{\cM_{g,n}(\mathbf{L})} \sum_{
    \begin{subarray}{1}
        \ Y \in \mathrm{Sub}_T(X), \\ \ 1 \leq |\chi(Y)| \leq A
    \end{subarray}} e^{-\lambda \ell_X(\p Y)} \ dX &\prec \sum_{m=1}^A m \cdot T^{4m+2} e^{(1-2\lambda)_+ T} \frac{V_{g,n}(\mathbf{L})}{g^m} \\ & \prec \frac{T^{4A + 2}e^{(1-2\lambda)_+ T} V_{g,n}(\mathbf{L})}{g}. \end{split} \eeq
    This completes the proof of Proposition \ref{Bounding exponential length function}.
    \ep

    \subsection{Concluding the proof of Theorem \ref{main theorem nonsimple curves} } The present subsection is devoted to the proof of Theorem \ref{main theorem nonsimple curves}; we begin with the following two lemmas.
    
    \bl \label{lemma 5.1} Fix $k < 0, c_T > 0$ and let $T = c_T \log g$. Then for any $n \geq 0$ and $\bf L \in \R_{\geq 0}^n$, there exists a constant $C = C(n,k)$ such that for any $g \geq 2$,
    \beq \EgnL \Bigg[\sum_{m=1}^\infty \sum_{\substack{\gamma \in \mathcal{P}^{ns}(X), \\ \ell_\gamma(X) \geq T}} \frac{e^{k \cdot m\ell_\gamma(X)}}{m(e^{m \ell_\gamma(X)} - 1)} \Bigg] \leq C(n,k) \cdot g^{1 - |k| c_T}. \eeq
    \el

    \bp Every $\gamma \in \cP^{ns}(X)$ satisfies $\ell_\gamma(X)> 4\sinh^{-1}(1)$ (see e.g. \cite[Theorem 4.2.2]{buser2010geometry}). Hence, there is a universal constant $c > 0$ such that for any $\gamma \in \cP^{ns}(X)$, we have
    \be \label{equivalence relation} c^{-1} \cdot \frac{e^{k \cdot \ell_\gamma(X)}}{e^{\ell_\gamma(X)} - 1}   \leq \sum_{m=1}^\infty \frac{e^{k \cdot m\ell_\gamma(X)}}{m(e^{m \ell_\gamma(X)} - 1)} \leq c \cdot \frac{e^{k \cdot \ell_\gamma(X)}}{e^{\ell_\gamma(X)} - 1}. \ene
    By Lemma \ref{Classical bound of geodesics}, for any $X \in \cM_{g,n}(\mathbf{L})$ and $L \geq 1$,
    \beq \#\{\gamma \in \cP^{ns}(X): \ell_\gamma(X) \leq L\} \leq \frac{2g+n-2}{2} e^{L+6}. \eeq
    It follows that for every $X \in \cM_{g,n}(\bf L)$,
    \beq \begin{split} \sum_{\substack{\gamma \in \mathcal{P}^{ns}(X), \\ \ell_\gamma(X) \geq T}} \frac{e^{k \cdot \ell_\gamma(X)}}{e^{\ell_\gamma(X)} - 1} &\leq \sum_{m = [T]}^\infty \frac{e^{km}}{e^{m} - 1} \cdot \#\{\gamma \in \cP_X^{ns}: \ell_\gamma(X) \leq m+1\} 
    \\ &\leq \sum_{m=[T]}^\infty \frac{e^{km}}{e^m - 1} \cdot \frac{2g+n-2}{2} e^{m+7} 
    \\ & \leq \frac{e^7}{1 - e^k} \cdot (2g+n-2) e^{k[T]} \\ &\leq \frac{e^{k + 7}}{1 - e^k} \cdot (2g+n-2) g^{k  c_T} \leq C(n,k) g^{1 - |k| c_T}, \end{split} \eeq
    where the penultimate inequality uses the definition $T = c_T \log g$. Taking the expectation of the above estimate over $\cM_{g,n}(\mathbf{L})$ and combining it with \eqref{equivalence relation}, we conclude 
    \beq \EgnL \Big[\sum_{m=1}^\infty \sum_{\substack{\gamma \in \mathcal{P}^{ns}(X), \\ \ell_\gamma(X) \geq T}} \frac{e^{k \cdot m\ell_\gamma(X)}}{m(e^{m \ell_\gamma(X)} - 1)} \Big] \leq C(n,k) \cdot g^{1 - |k| c_T}. \eeq
    This completes the proof.
    \ep

    \bl \label{lemma 5.2} 
    
    Fix $n\geq 0$ and suppose \(|\mathbf L|^2=o(g)\) as \(g\to\infty\). Let $c_T > 0$ be fixed and define $T = c_T \log g$. Then for any $\eps > 0$, 
    \beq \EgnL \Big[\sum_{m=1}^\infty \sum_{\substack{\gamma \in \mathcal{P}^{ns}(X), \\ \ell_\gamma(X) \leq T}} \frac{1}{m(e^{m \ell_\gamma(X)} - 1)} \Big] \prec \frac{(c_T \log g)^2}{g}  + \frac{(c_T \log g)^{16 c_T + 3}}{ g^{1 - c_T \eps} }. \eeq
    
    \el

    \bp Since every $\gamma \in \cP^{ns}(X)$ satisfies $\ell_\gamma(X) > 4 \sinh^{-1}(1)$, there is a universal constant $c > 0$ such that for any $\gamma \in \cP^{ns}(X)$, we have
    \be \label{equivalence relation 2} c^{-1} \cdot e^{- \ell_\gamma(X)} \leq \sum_{m=1}^\infty \frac{1}{m(e^{m \ell_\gamma(X)} - 1)} \leq c \cdot e^{- \ell_\gamma(X)}. \ene
    For each non-simple closed geodesic, write $Y$ for its associated filling subsurface introduced in Section \ref{filling construction}. We then rewrite the sum as
    \beq \sum_{\substack{\gamma \in \mathcal{P}^{ns}(X), \\ \ell_\gamma(X) \leq T}} e^{-\ell_\gamma(X)} = 2\sum_{Y \in \mathrm{Sub}_{T}(X)} \sum_{\gamma \text{ fills } Y} e^{- \ell_\gamma(X)} \mathbf{1}_{\{\ell_\gamma(X) \leq T\}}. \eeq
    We divide the summation on the right-hand side into two disjoint subfamilies, distinguished by the bound $|\chi(Y)| < A$ versus $|\chi(Y)| \geq A + 1$. Since $\gamma$ fills in $Y$ and $\ell_\gamma(X) \leq T$, we have $\ell(\p Y) \leq 2 \ell_\gamma(Y)$ and $\area(Y) \leq 4T$. Combining these two estimates with Lemma \ref{Classical bound of geodesics} gives
    \beq \begin{split} \sum_{
    \begin{subarray}{1}
        Y \in \mathrm{Sub}_T(X), \\  |\chi(Y)| \geq A+1
    \end{subarray}} \sum_{\gamma \text{ fills } Y } e^{- \ell_\gamma(Y)} \mathbf{1}_{\{\ell_\gamma(X) \leq T\}} & \leq \sum_{
    \begin{subarray}{1}
        Y \in \mathrm{Sub}_T(X), \\  |\chi(Y)| \geq A+1
    \end{subarray}} e^{-\frac{\ell(\p Y)}{2}} \#_f(Y,T) 
    \\ & \leq \frac{e^6}{\pi} T e^T \sum_{
    \begin{subarray}{1}
        Y \in \mathrm{Sub}_T(X), \\  |\chi(Y)| \geq A+1
    \end{subarray}} e^{- \frac{\ell(\p Y)}{2}}. \end{split} \eeq
    We next handle subsurfaces with relatively small Euler characteristics \(1\le|\chi(Y)|\le A\):
    \be \begin{split} & \quad \ \sum_{
    \begin{subarray}{1}
        \ Y \in \mathrm{Sub}_T(X), \\ \ 1 \leq |\chi(Y)| \leq A
    \end{subarray}}  \sum_{\gamma \text{ fills } Y } e^{-\ell_\gamma(X)}\mathbf{1}_{\{\ell_\gamma(X) \leq T\}} 
    \\ & \leq \sum_{m=1}^{[T]+1} \sum_{
    \begin{subarray}{1}
        \ Y \in \mathrm{Sub}_{m}(X), \\ \ 1 \leq |\chi(Y)| \leq A
    \end{subarray}} \sum_{\gamma \text{ fills } Y } e^{-\ell_\gamma(X)}\mathbf{1}_{\{m-1 < \ell_\gamma(X) \leq m\}}
    \\ & \leq \sum_{m=1}^{[T]+1} \sum_{
    \begin{subarray}{1}
        \ Y \in \mathrm{Sub}_{m}(X), \\ \ 1 \leq |\chi(Y)| \leq A
    \end{subarray}} e^{-(m-1)} \cdot \#_f(Y,m). \end{split} \label{small Euler characteristics final bound} \ene
    By Lemma \ref{main counting WX22}, for any $0 < \eps < 1/2$ and subsurface $Y \subset X\in \cM_{g,n}(\mathbf {L})$ with $ |\chi(Y)| \leq A$, there is a constant $c(\eps,A)>0$ depending only on $\eps, A$ such that for any $m \geq 1$,
    \beq \#_f(Y,m) \leq c(\eps,A) e^m e^{- \frac{(1-\eps) \ell(\p Y)}{2}}. \eeq
    Substituting this bound into \eqref{small Euler characteristics final bound} and performing the summation over $m$ yields
    \beq \sum_{
    \begin{subarray}{1}
        \ Y \in \mathrm{Sub}_T(X), \\ \ 1 \leq |\chi(Y)| \leq A
    \end{subarray}} \sum_{\gamma \text{ fills } Y } e^{-\ell_\gamma(X)}\mathbf{1}_{\{\ell_\gamma(X) \leq T\}}  \leq c(\eps,A) T \sum_{
    \begin{subarray}{1}
        \ Y \in \mathrm{Sub}_T(X), \\ \ 1 \leq |\chi(Y)| \leq A
    \end{subarray}} e^{-\frac{(1 - \eps)\ell(\p Y)}{2}}. \eeq
    Combining all these inequalities and taking the expectation, we have
    \beq \begin{split} \EgnL \Big[\sum_{\substack{\gamma \in \mathcal{P}^{ns}(X), \\ \ell_\gamma(X) \leq T}} e^{-\ell_\gamma(X)} \Big] 
    & \prec \EgnL \Big[ T e^T \sum_{ \begin{subarray}{1}
        Y \in \mathrm{Sub}_T(X), \\  |\chi(Y)| \geq A+1
    \end{subarray}} e^{- \frac{\ell(\p Y)}{2}} \Big]
    \\ & +  \EgnL \Big[T 
    \sum_{
    \begin{subarray}{1}
        \ Y \in \mathrm{Sub}_T(X), \\ \ 1 \leq |\chi(Y)| \leq A
    \end{subarray}} e^{-\frac{(1-\eps)\ell(\p Y)}{2}} \Big]. \end{split} \eeq
    Applying Proposition \ref{Bounding exponential length function}, we deduce 
    \beq \EgnL \Big[\sum_{\substack{\gamma \in \mathcal{P}^{ns}(X), \\ \ell_\gamma(X) \leq T}} e^{- \ell_\gamma(X)} \Big] \prec \frac{T^2 e^{4T}}{g^{A+1}}  + \frac{T^{4A + 3} e^{\eps T}}{g}. \eeq
    Recall our choice \(T = c_T\log g\). Substitute \(A = 4c_T\) to obtain
    \beq \EgnL \Big[\sum_{\substack{\gamma \in \mathcal{P}^{ns}(X), \\ \ell_\gamma(X) \leq T}} e^{-\ell_\gamma(X)} \Big] \prec \frac{(c_T \log g)^2}{g}  + \frac{(c_T \log g)^{16 c_T + 3}}{ g^{1 - c_T \eps} }. \eeq
    This, together with the equivalence relation \eqref{equivalence relation 2}, completes the proof.
    \ep
    
    Now we are ready to prove Theorem \ref{main theorem nonsimple curves}.

    \bp[Proof of Theorem \ref{main theorem nonsimple curves}] Fix $\kappa > 0$, correspondingly $k < 0$. By \eqref{LW26-0} we know that
    \beq \mu_{\bar X}^\kappa(\cC^{ns}(\bar X)) = \sum_{m=1}^\infty \sum_{\gamma \in \cP^{ns}(X)} \frac{e^{k \cdot m \ell_\gamma(X)}}{m(e^{m \ell_\gamma(X)} - 1)}. \eeq
    Set $c_T = 2/|k|$ so that \(1-|k|c_T=-1\). By Lemma \ref{lemma 5.1},
    \be  \EgnL \Big[\sum_{m=1}^\infty \sum_{\substack{\gamma \in \mathcal{P}^{ns}(X), \\ \ell_\gamma(X) \geq T}} \frac{e^{k \cdot m\ell_\gamma(X)}}{m(e^{m \ell_\gamma(X)} - 1)} \Big] \prec \frac{1}{g}. \label{nonsimple endgame case1} \ene
    Moreover, Lemma \ref{lemma 5.2} guarantees that for any $0 < \eps < 1$,
    \be \EgnL \Big[\sum_{m=1}^\infty \sum_{\substack{\gamma \in \mathcal{P}^{ns}(X), \\ \ell_\gamma(X) \leq T}} \frac{1}{m(e^{m \ell_\gamma(X)} - 1)}  \Big] \prec \frac{1}{g^{1-\eps}}. \label{nonsimple endgame case2} \ene 
    Combining \eqref{nonsimple endgame case1} and \eqref{nonsimple endgame case2} yields
    \be \EgnL \Big[ \mu_{\bar X}^{\kappa}(\cC^{ns}(\bar X)) \Big] \prec \frac{1}{g^{1-\eps}}, \ \ \forall \ 0 < \eps < 1. \ene
    This completes the proof of Theorem \ref{main theorem nonsimple curves}.
    \ep

\bibliographystyle{amsalpha}
\bibliography{ref}

\end{document}